\documentclass{amsart}

\newtheorem{theorem}{Theorem}[section]
\newtheorem{lemma}[theorem]{Lemma}

\theoremstyle{definition}
\newtheorem{definition}[theorem]{Definition}
\newtheorem{example}[theorem]{Example}

\usepackage{amsthm,amssymb,amsmath,amsfonts}
\usepackage{graphicx}
\usepackage{subcaption}
\usepackage[export]{adjustbox}
\usepackage{wrapfig}
\usepackage{relsize,exscale}
\usepackage{algorithm}
\usepackage{program}

\theoremstyle{remark}
\newtheorem{remark}[theorem]{Remark}

\numberwithin{equation}{section}

\begin{document}

\title{The Limit Space of Self-similar Groups and Schreier graphs}

\author{bozorgmehr vaziri}
\author{farhad rahmati}
%\dedicatory{This paper is dedicated to soheila sarabandi}
\begin{abstract}

The present paper investigates the limit $G$-space $\mathcal{J}_{G}$ generated by
the self-similar action of automatic groups on a regular rooted tree. The
limit space $\mathcal{J}_{G}$ is the Gromov-Hausdorff limit of the family of Schreier
graphs $\Gamma_{n}$; therefore, $\mathcal{J}_{G}$ can be approximated by Schreier graphs on
level $n$-th when $n$ tends to infinity. We propose a computer program
whose code is written in Wolfram language computes the adjacency matrix of Schreier graph $\Gamma_{n}$ at each specified level of the
regular rooted tree. In this paper, the Schreier graphs corresponding to each automatic group is computed by applying the program to some
collection of automata groups, including classic automatic groups.
\end{abstract}

\maketitle

\section{Introduction}\label{sec1}

The idea of the boundary at infinity is due to Gromov, who attached to a proper hyperbolic metric space $(X,d)$ a compact metrizable space $\partial X$ whose points are classes of asymptotic geodesic rays in $X$,  reflecting all coarse geometry (i.e., quasi-isometric invariants) of the metric space $X$. This attachment has played a central role in recent advances in geometry, geometric group theory, ergodic theory, and dynamical systems called \textit{limit space} (\textit{boundary} or \textit{boundary at infinity}). When a discrete group $G$ acts geometrically (co-compactly properly discontinuous by isometry) on the metric space $X$, besides the geometric structure of the space $X$, the boundary $\partial X$ also reflects the asypmtotic structure of the group $G$ (like the growth rate, the complexity of the word problem and the Hausdorff
dimension of the limit space),
 since the $G$-action $G\curvearrowright X$ extends to a homeomorphic $G$-action on the boundary $G\curvearrowright \partial X$. When $X$ is a complete $n$-dimensional Riemannian manifold of negative sectional curvature, then the boundary $\partial X$ is homeomorphic to the $(n-1)$-sphere, such that for given a base point $x_{0} \in X$, there is a homeomorphism by the map $T_{x_{0}}X\rightarrow \partial X$ which associates to each tangent unit vector $u$ at the base point the class of the geodesic ray $c:[0,\infty)\rightarrow X$  issues from the base point with velocity vector $u$. Furthermore, there is a $G$-invariant measure $\mu$ on the boundary called harmonic measure, making $\partial X$ a dynamical system $(\partial X, G, \mu)$.

The limit space associated with a negatively curved metric space possesses impressive rich geometry,
topology, and dynamics, which rised the attractions from different realms of mathematics and modern physics in the last two decades. For instance, the author \cite{yue1996ergodic} uses ergodic theory at infinity of an arbitrary discrete isometry group $\Gamma$ acting on any Hadamard manifold $\mathbb{H}$ of pinched variable negative curvature to measure the thickness of the galaxy. In addition, the dynamic of the geodesic flow on a negatively curved manifold consider an ideal chaotic behavior in quantum dynamical systems.
The limit space also is used to model the long-time behavior of the chaotic physical systems to estimate the thermodynamic equilibrium (see \cite{bourgain2003entropy},\cite{anantharaman2007entropy},\cite{jezouin2013quantum}). 
\\
The notation of limit space (or asypmtotic boundary) in geometry group theory faced a significant advance after impressive works of Gromove \cite{gromov1980hyperbolic}\cite{gromov1981groups}\cite{gromov1992asymptotic}, also
since 1990 the theory of automatic groups acting on a homogeneouse rooted tree and their limit space extensively developed by Grigorchuk, Bartholdi, Nekrashevish, Sidiki, Bondarenko (\cite{bartholdi2003fractal}\cite{grigorchuk2011some}\cite{grigorchuk2014self}\cite{grigorchuk2006asymptotic}\cite{nekrashevych2005self}
\cite{bartholdi2006automata}\cite{bondarenko2011ends})
To each self-similar contracting group,  \cite{nekrashevych2005self} associates a limit-space is defined as a quotient of the space $A^{-\omega}$ of the left-infinite words over a nite alphabet $A$, modulo a $G$-invariant asymptotic equivalence relation, that can describe using the Moore diagram of the automaton generating the group. This limit space is approximate by a sequence of finite graphs, namely the finite Schreier graphs associated with the automorphism action of the automatic group $G$ on the regular rooted tree $A^{*}$.
It is shown when $G$ acts by contracting on $A^{*}$, then $G$ is a word hyperbolic group by the word metric, and there is an algorithm to decide whether two automatic fractal groups, $G_{1}$ and $G_{2}$, are isomorphic subgroups in the automorphism group $Aut(A^{*})$. Accordingly, the contracting condition prepares an algorithmic way to make the algebraic structure of automatic groups computable. These additional structures (self-similarity and contracting) provide a significant advantage; it paves the way for bringing the automatic groups and associated limit space $\mathcal{J}_{G}$ into computer science and computational geometry. Moreover \cite{nekrashevych2005self} proves the limit space $\mathcal{J}_{G}$ has finite topological dimension (i.e., complexity) if and only if $G$ has a finite non-empty nucleus (contracting).

Authors \cite{bartholdi2006automata} explore the connections between automata, groups, limit spaces of self-similar actions, and
tilings. In particular, they show how a group acting “nicely” on a tree gives rise to a self-covering of
a topological groupoid and how the group can be reconstructed from the groupoid and its covering and
The connection is via finite-state automata. V. Nekrashevych and A. Teplyaev \cite{nekrashevych2008groups} investigate the relation between the Analysis of fractals and self-similar Groups. The random walks on a crystal or quasi-crystal are formed by the direct limit of the random walks on a sequence of metric spaces studied by \cite{kigami2001analysis}\cite{lindstrom1990brownian}.

The current research considers the limit space $\mathcal{J}_{G}$, where $G$ is a contracting fractal group, and tries to supply a computational mechanism in the form of a computer program to approximate $\mathcal{J}_{G}$ by its corresponding Schreier graphs. In this direction, the package AutomGrp in the GAP program environment is created by Muntyan and Savchuk \cite{AutomGrp1.3.2}. It provides methods or computations with groups and semigroups generated by finite automata or given by wreath recursion, as well as with their finitely generated subgroups and elements. The other is a computer program presented in the Ph.D. thesis \cite{perez2020structural} written by Perez in 2019 that produces Schreier graphs of automatic groups.

\begin{figure}[h!]
\centering
  \begin{subfigure}[b]{0.6\linewidth}
\centering
    \includegraphics[width=1.1\textwidth5]{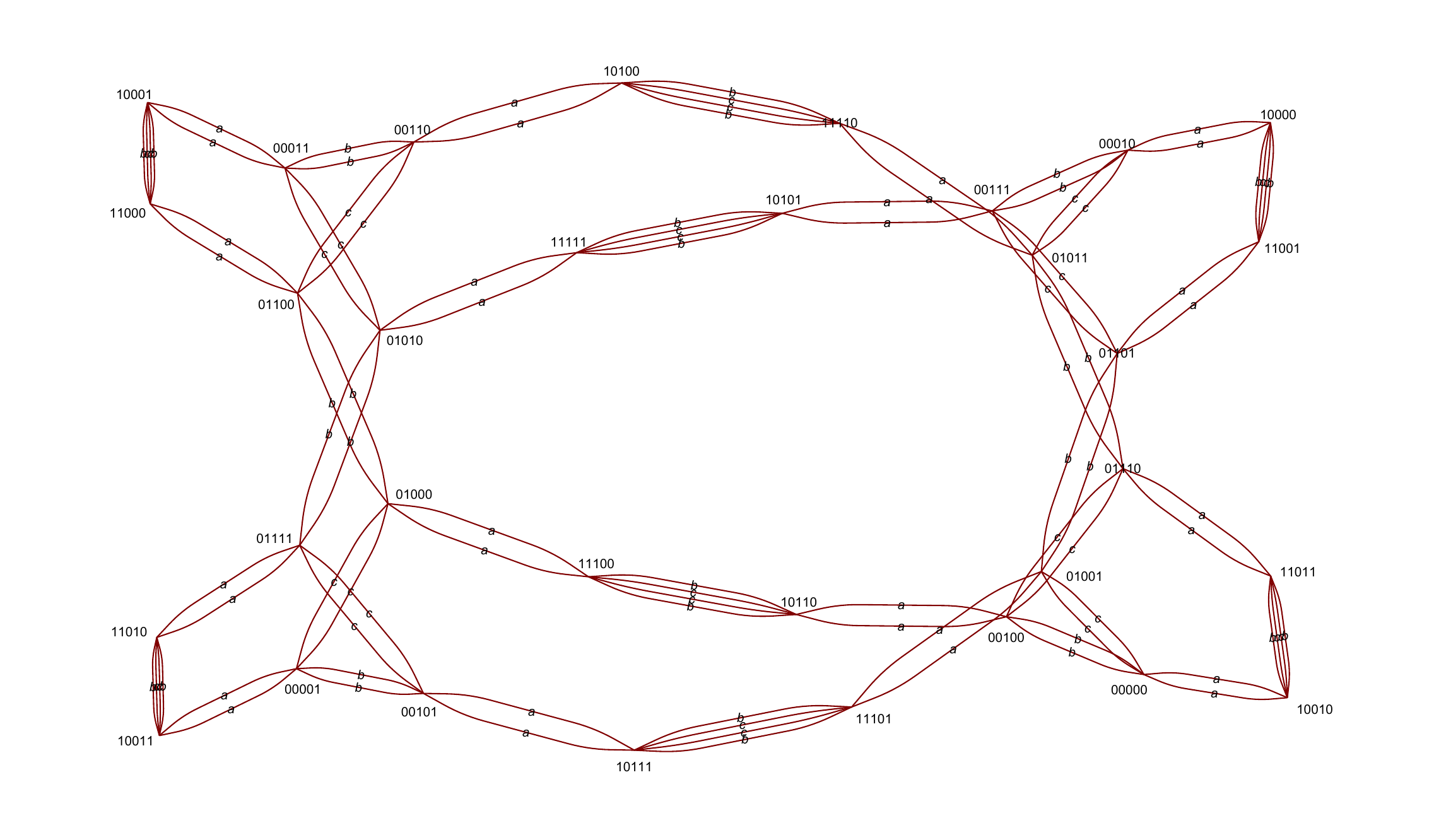}
    \caption{Schreier graph level 5}
  \end{subfigure}
  \begin{subfigure}[b]{0.45\linewidth}
\centering
    \includegraphics[width=0.9\textwidth5]{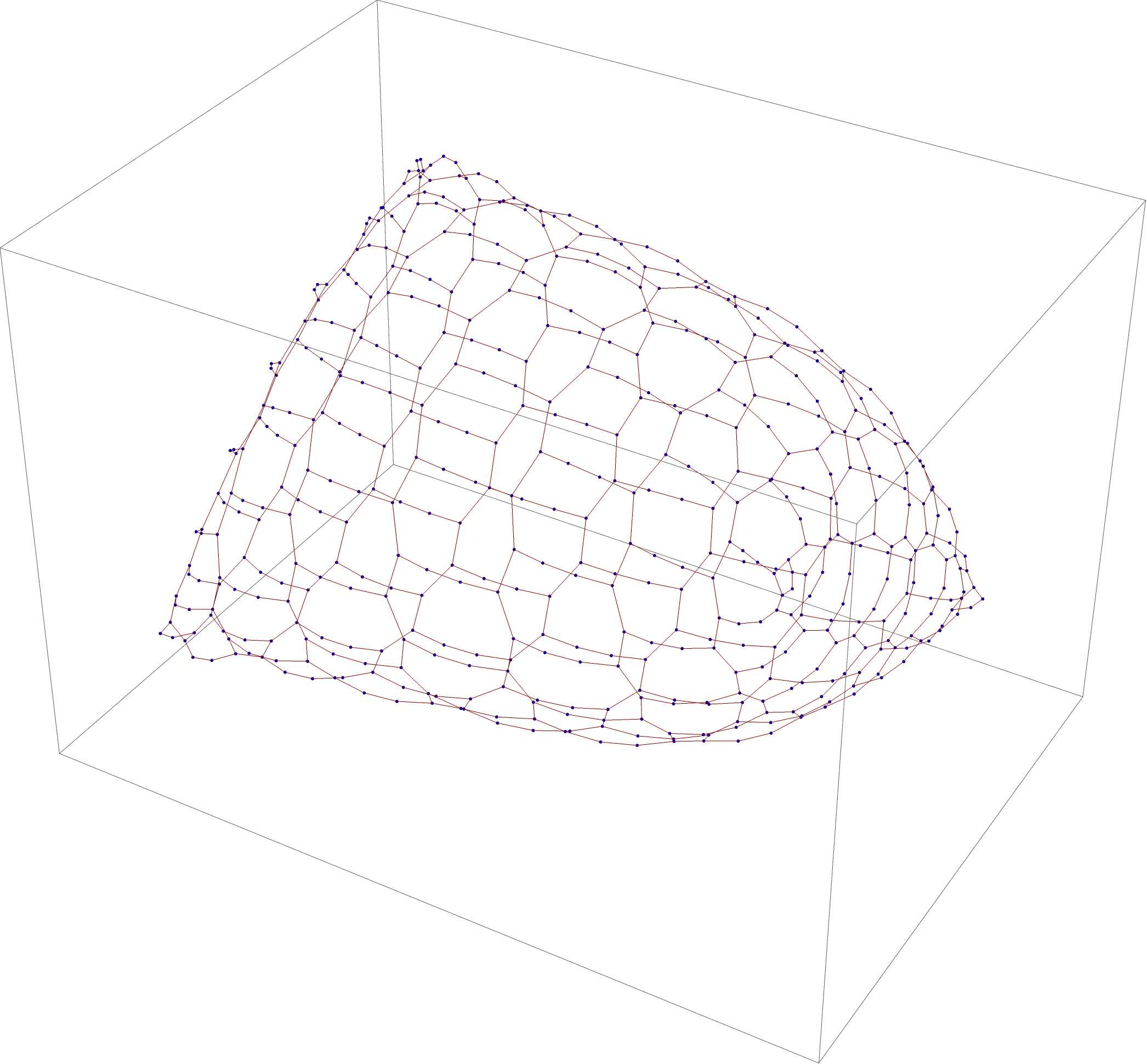}
    \caption{simplicial Schreier graph level 9}
  \end{subfigure}
\begin{subfigure}[b]{0.45\linewidth}
\centering
    \includegraphics[width=0.9\textwidth]{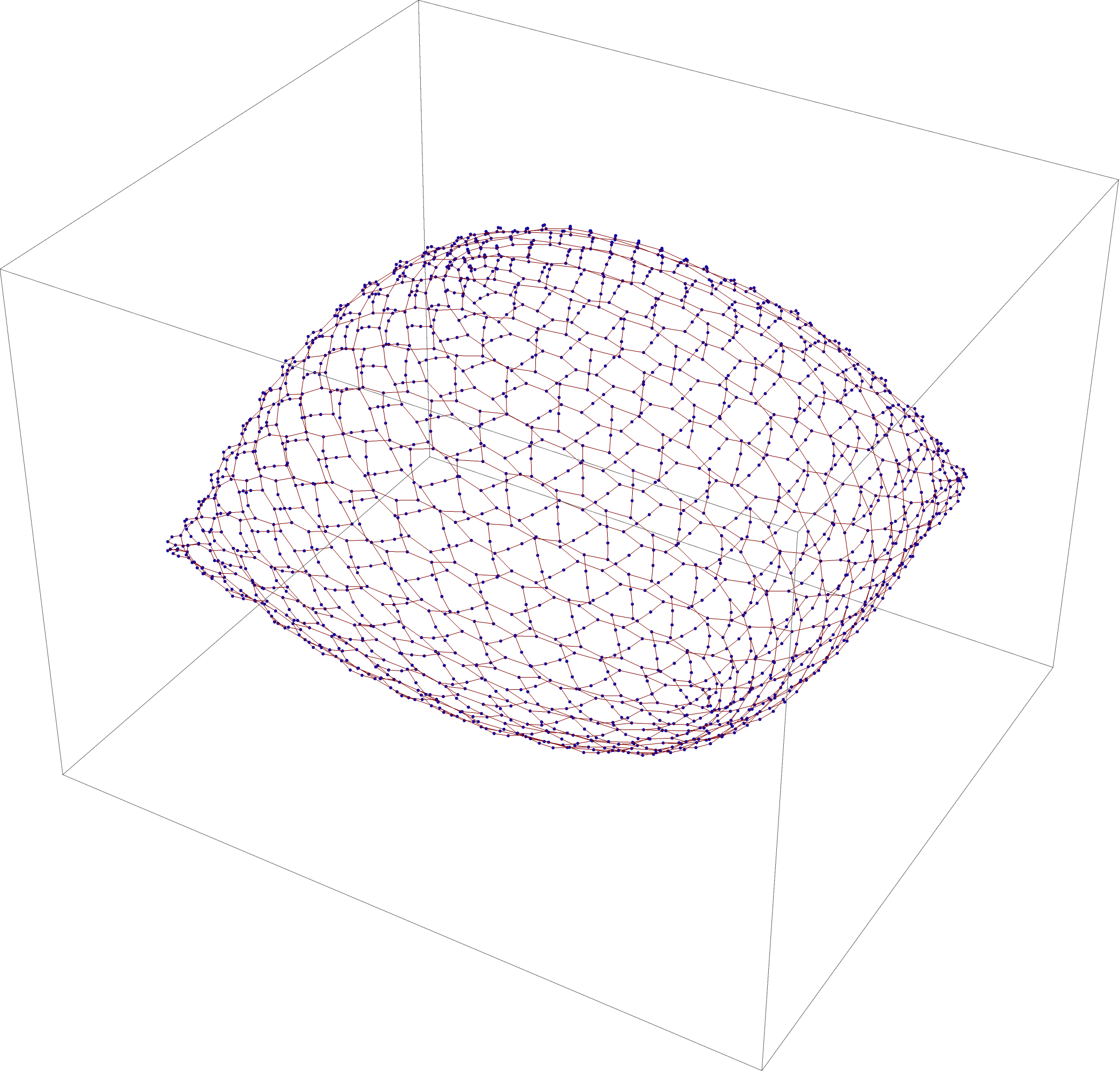}
     \caption{simplicial Schreier graph level 11}
  \end{subfigure}
  \caption{The Schreier graphs generated by automata $a=\sigma(b,b) , b=(c,a) , c=(a,a)$ is a virtually $\mathbb{Z}^{3}$ group}
  \label{virtually Z3}
\end{figure}

The limit dynamical system $(\mathcal{J}_{G}, s ,\mu)$ is naturally used to model chaotic systems with self-similar structures; that is why the limit-space is highly capable of applying to different areas of science. Before anything, we require to approximate the incidence structure of $\mathcal{J}_{G}$. Fortunately, the Schreier graphs of a self-similar action $(G, A)$ make this possible. In other words, Schreier graphs provide an algorithmic way to approximate the limit space by an arrangement of finite graphs (see figure \ref{virtually Z3}).

We consider the demand for computational tools and propose a procedure for calculating the adjacency structure of Schreier graphs for the arbitrary $n$-th level of the rooted tree. This procedure is executed by a computer program whose codes (given in the last section \ref{Program codes}) are composed in Wolfram Mathematica language. Intuitively, the program is a machine that receives the information of finite states automata consisting of The adjacency structure of the Moor diagram associated with the automaton and the permutation action of each state $q \in Q$ of automata on the alphabet set $A$. Then the machine creates the adjacency matrix of the Schreier graph $\Gamma_{n}$ for an arbitrary integer $n>0$.

The boundary space $\mathcal{J}_{G}$ can also be defined as an inverse limit of the
sequence of covering projection $\Gamma_{n}\rightarrow \Gamma_{n-1}$, where $\Gamma_{n}$ is a ball of radius $n$ around the identity element in the Cayley graph of the group. Moreover, the sequence converges to boundary space by Gromov-Hausdorff distance (see \cite{bridson2013metric}).

As can be seen, there are various approaches to dealing with
the limit-space of negatively curved metric spaces since the limit space possesses different structures that are well behaved and compatible with each other. For instance, the boundary has inverse limit topology (pro-discrete topology). The visual metrics $d_{x_{0}}^{\alpha}$ for real parameter $\alpha$  and base point $x_{0}$ yield canonical quasi-conformal and conformal structures, and a $G$-invariant measure also admits ergodic dynamical systems and canonical orbispace structures.
\begin{figure}[h!]
  \centering
  \includegraphics[width=0.9\linewidth]{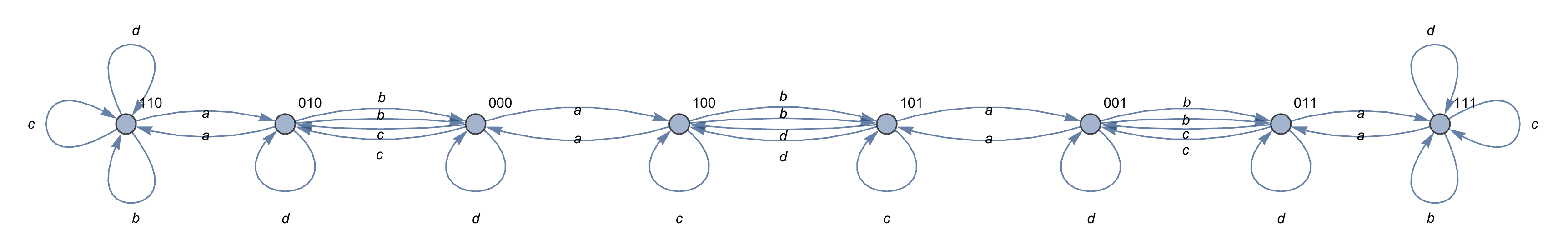}
 
  \caption{The Schreier graph of the first Grigorchuk group at level 3-th }
  \label{The Schreier graphs by first Grigorchuk group}
\end{figure}

The asymptotic equivalence relation on the boundary $\mathcal{X}_{G}:= A^{-\omega}$ is determined by the infinite directed paths in the Moor diagram of automata. 
It is not hard to see that the contracting assumption yields a $G$-invariant closed asymptotic equivalence relation $\sim$ on the boundary $\mathcal{X}_{G}$, making the limit-space $\mathcal{J}_{G}$ a compact Hausdorff quotient space $\mathcal{X}_{G}=A^{-\omega}/\sim$ with finite topological dimension $\leq \vert \mathcal{N}\vert - 1$ and the embedding
\begin{center}
$\mathcal{J}_{G}\hookrightarrow \mathbb{R}^{2\vert \mathcal{N}\vert -1}$
\end{center}
Nekrashevich \cite{nekrashevych2005self} proves that when $G$ is a contractive and self-replacing automatic group, $G$ is generated by the finite nucleus $\mathcal{N}$, and the limit space $\mathcal{J}_{G}$ is a path-connected and locally path-connected compact metric space accompanied by the $G$-bundle $\mathcal{X}_{G}$ and projecting map $\pi :\mathcal{X}_{G}\rightarrow\mathcal{J}_{G}$, that provides a maximal orbifold structure on the limit space.
Indeed, the self-similar structure on limit spaces is produced by extended self-similar action of the countable automatic group $G $ on the spherically homogenous rooted tree $A^{*}$.

The geometric action of $G$ on the \textit{self-similarity graph} $\sum(G, S, X)$ (see definition), which is a proper hyperbolic space demonstrated by \cite{nekrashevych2005self}, implies that the group $G$
to be a word hyperbolic group in the sense of Gromov \cite{gromov1987hyperbolic}. 
Consequently, the limit space $\mathcal{J}_{G}$ can interpret as the boundary of a proper hyperbolic space $\sum(G, S, X)$ with homeomorphism $\partial G\rightarrow\partial\sum(G,S,X)$ induced by the quasi-isometry between $G$ and $\sum(G, S,X)$. For more detail about boundary of hyperbolic groups see \cite{kapovich2002boundaries}

\begin{example}\cite{d2009schreier}
The Basilica group is a torsion free weakly branch group G without free subgroups defined by was introduced by R. Grigorchuk and A. Zuk in \cite{grigorchuk2002torsion}, ˙
where they show that it does not belong to the closure of the set of groups of subexponential growth
under the operations of group extension and direct limit. L. Bartholdi \cite{bartholdi2010amenability} showed it
to be amenable, making Basilica the first example of an amenable but not subexponentially amenable
group, and there exists therefore a natural way to associate to it a
compact limit space homeomorphic to the well-known Basilica fractal (see \ref{Schreier graph of Basilica}).
Let $G$ be group generated by finite automata with contracting and fractal action on binary tree
\begin{center}
$a=\sigma(b , id)  ,  b=(a , id)  ,  id=(id ,id)$   
\end{center}

It is shown that the \textit{iterated monodromy group} $IMG(z^{2}-1)$ which it's limit $G$-space is homeomorphic to the Julia set of polynomial $z^{2}-1$\cite{nekrashevych2005self}.

 \begin{figure}
 \centering
    \includegraphics[width=0.9\linewidth]{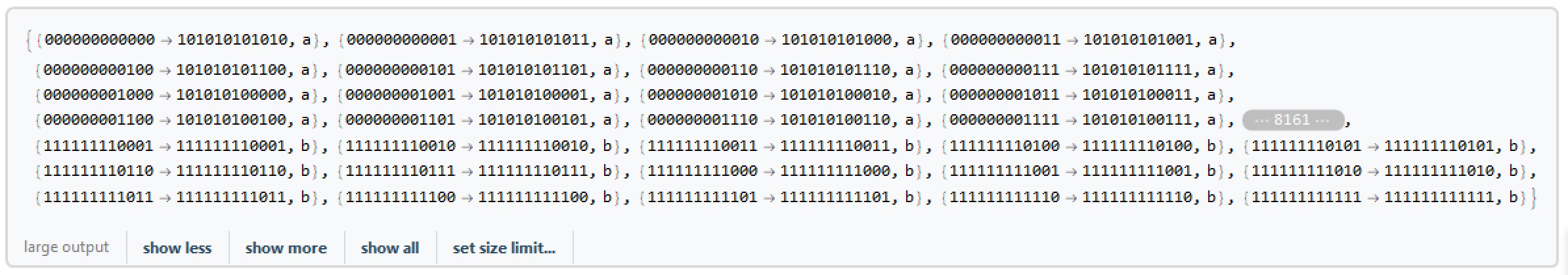}
     \caption{ The calculation of the adjacency structure of the Basilica group
 $IMG(z^{2}-1)$ at level 12}
  \label{computation of Basilica group}
\end{figure}

 \begin{figure}
 \centering
    \includegraphics[width=0.8\linewidth]{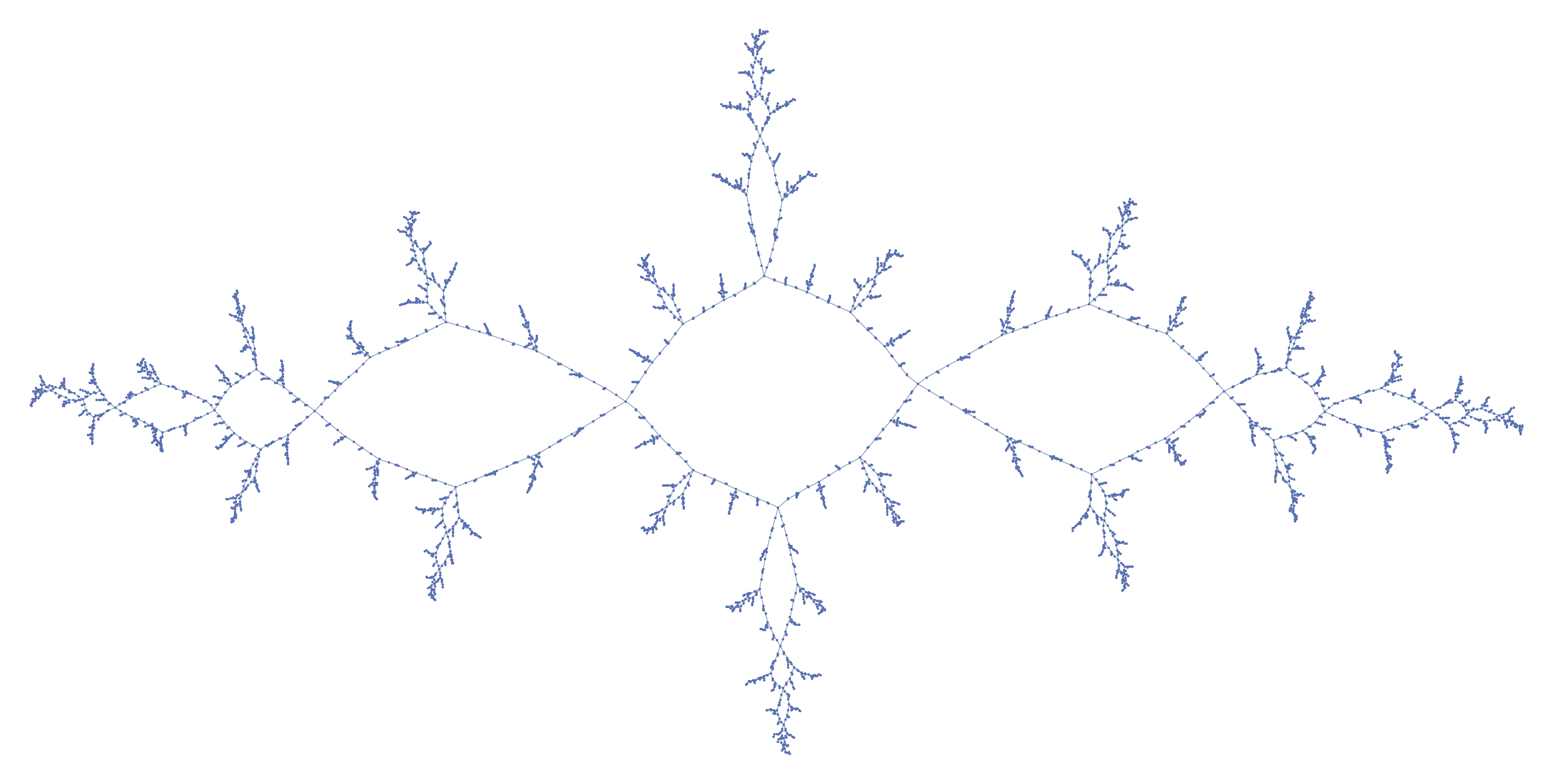}
     \caption{The Schreier graph of Basilica group $IMG(z^{2}-1)$ at level 12}
  \label{Schreier graph of Basilica}
\end{figure}
\end{example}

  \begin{figure}[h]
  \begin{subfigure}[b]{0.45\linewidth}
  \centering
    \includegraphics[width=1\textwidth]{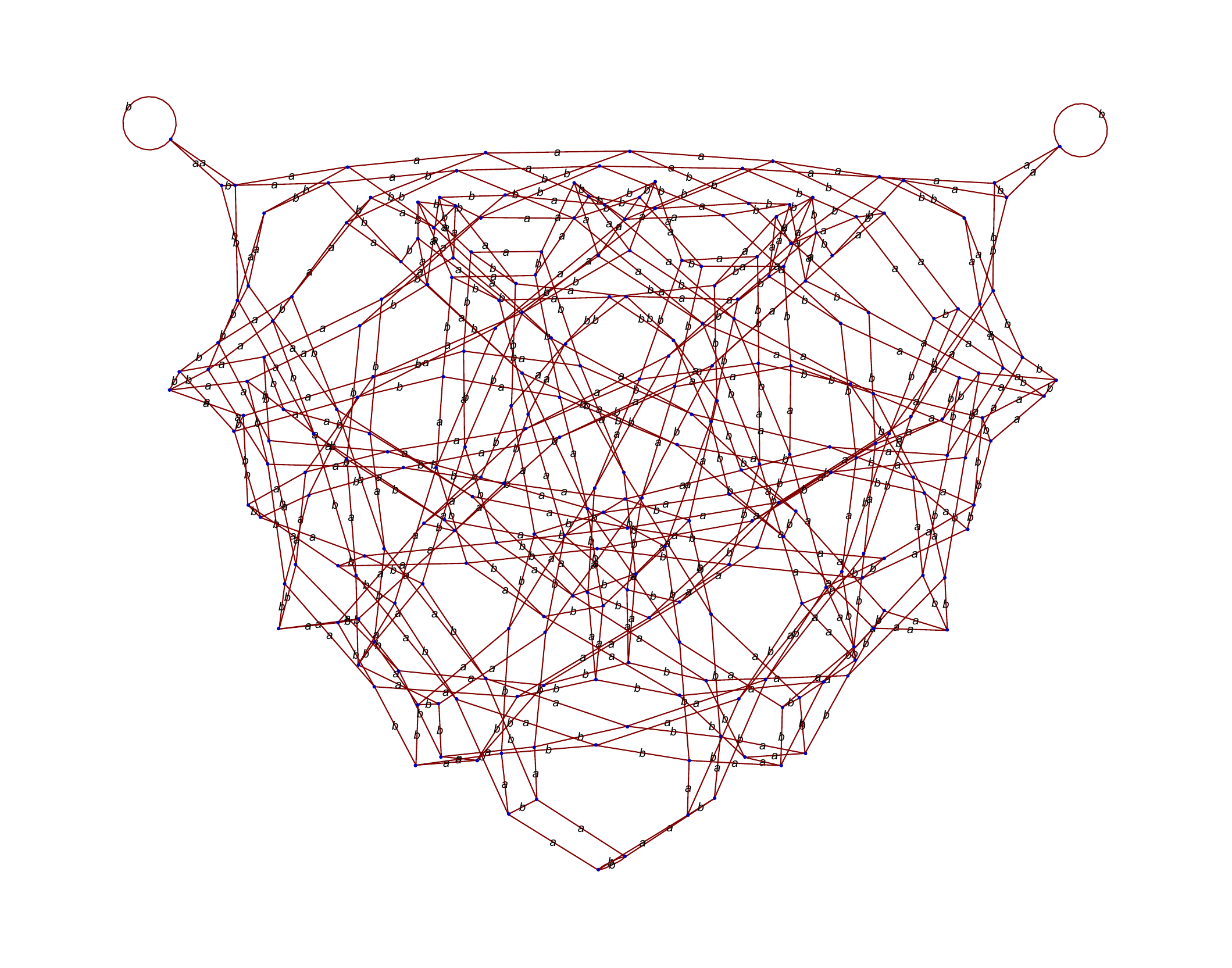}
    \caption{
level 8 }
  \end{subfigure}
\begin{subfigure}[b]{0.45\linewidth}\centering
    \includegraphics[width=1\textwidth]{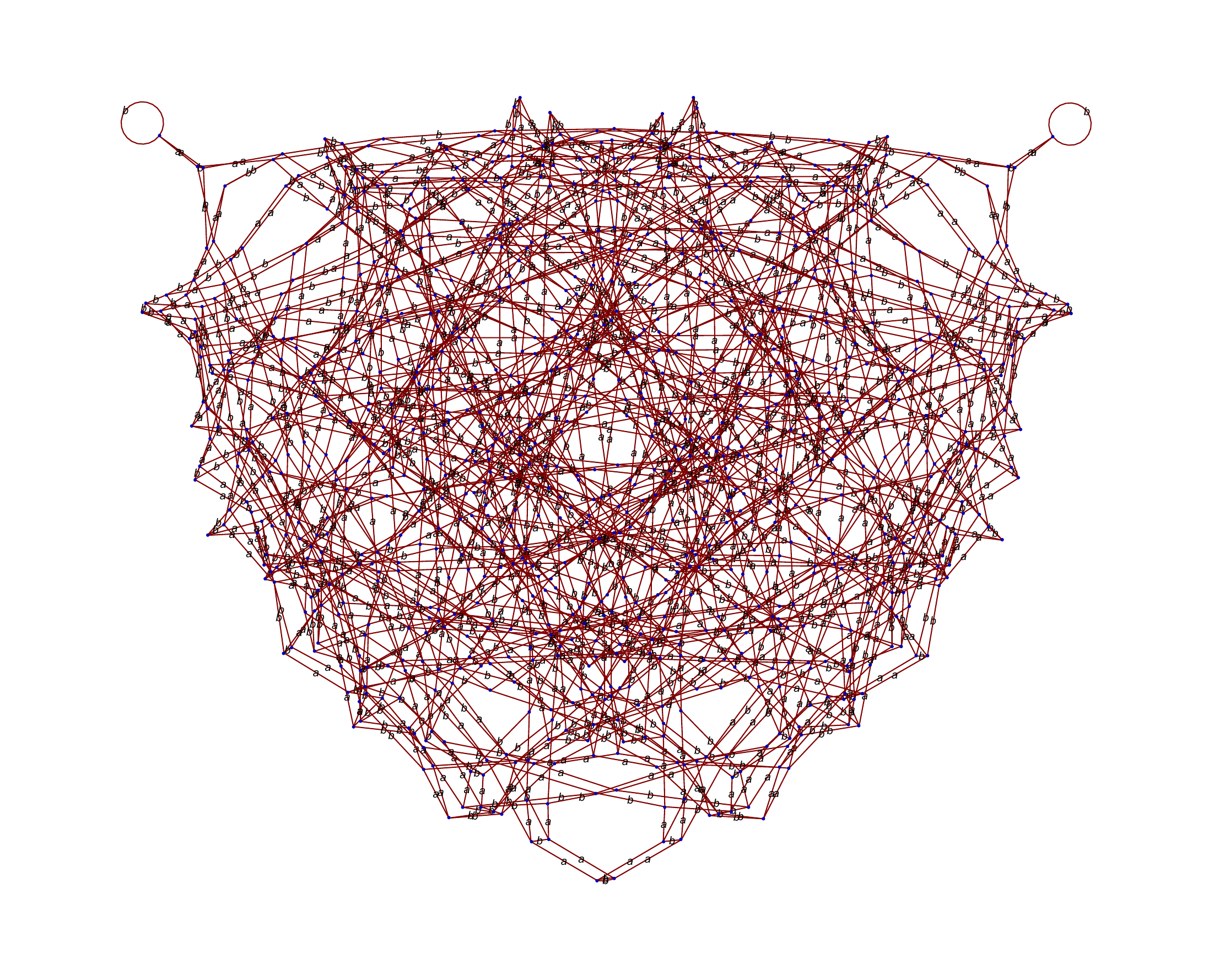}
     \caption{
level 10  }
  \end{subfigure}
  \caption
{The Schreier graph of the Lamplighter group $\mathbb{Z}_{2}\wr\mathbb{Z}$ generated with the wreath recursion $a=\sigma (b,a) , b=(b,a)$.}
  \label{The Schreier graph of Lamplighter group}
\end{figure}
%\begin{figure}[h!]
%  \centering
%  \includegraphics[width=0.7\linewidth]{automataLamplightergrouplevel8.pdf}
% \caption{The tail graph of Lamplighter group $\mathbb{Z}_{2}\wr\mathbb{Z}$ at level 8}
%  \label{The tail graph of Lamplighter group}
%\end{figure}

The self-similar contracting groups appear naturally in studying expanding (partial) self-coverings of topological spaces and orbispaces as their iterated monodromy groups. When the group is the iterated monodromy group of an expanding partial self covering map $f$, \cite{nekrashevych2005self} shows the limit space is homeomorphic to the Julia set $J(f)$. Finite Schreier graphs corresponding to the group $G$,  form a sequence of combinatorial approximations to the limit space. If the open conditions satisfy, the limit-space is homeomorphic to a post-critically finite self-similar fractal, being generated by bounded automaton guarantees that the limit space is post-critically finite. Authors \cite{bartholdi2006automata} explore the connections between automata, groups, limit spaces of self-similar actions, and tilings. In particular, they show how a group acting nicely on a tree gives rise to a self-covering of a topological groupoid and how the group can be reconstructed from the groupoid and its covering via finite-state automata. V. Nekrashevych and A. Teplyaev \cite{nekrashevych2008groups} investigate the relationship between the Analysis of fractals and self-similar Groups. The random walks on a crystal or quasi-crystal are formed by the direct limit of the random walks on a sequence of metric spaces studied by \cite{kigami2001analysis}\cite{lindstrom1990brownian}
\\
The structure of the present paper is as follows: 
The introduction\ref{sec1} briefly describes the theory of contracting self-similar groups acting on a spherically homogeneous rooted tree and their limit G-space. 
In the second section\ref{preliminary}, the definitions and known results about the Schreier graphs and limit-space are given.
In section \ref{result}, the program is applied to the collection of automatons and explores the incidence structure of the family of Schreier graphs.
The collection contains the first Grigorchuk group, Basilica group, the free group $\mathbb{F}_{3}$ of rank 3,Sierpinski group, the lamplighter group mother groups, longe range group,Hanoi tower group, and some automatic groups that are expressed in \cite{grigorchuk2014self}\cite{bondarenko2008classification}.
In the final section, the code written in Wolfrom Mathematica language is presented, which is used to construct the Schreier graphs and calculate their adjacency structure.
\section{Preliminary}\label{preliminary}
This section reviews the definitions and theorems related to the self-similar groups generated by closed finite-state automata, their Schreier graphs, and the limit space.We have left theorems and propositions without proof, guiding the reader to the references mentioned.
The main resources that are used in this section are \cite{bartholdi2003fractal}\cite{bartholdi2006automata}\cite{nekrashevych2005self}\cite{nekrashevych2008groups}.

%\begin{figure}[h!]
%  \centering
%   \includegraphics[width=1\linewidth]{com2853.pdf}
%  \caption{Schreier graph of $IMG(\dfrac{z-1}{z+1}^{2})\cong$ group generated by the wreath recurssion $a=\sigma(c,c) , b=\sigma (b,a) ,c=(c,c)$.}
%  \label{The spectrum of Laplacian operator on the Schreier graphs by the Lamplighter group on complex plane}
%\end{figure}

\subsubsection{Group generated by automaton}
Let $A$ is a alphabet set, the free monoid $A^{*}$ with empty set consider as root tree and $A^{\omega}$ denote the \textit{end} (\textit{boundary}) of rooted tree consist of geodesic paths issuing from a point  $x\in A^{*}$ in the tree.
\begin{definition}
 The \textit{automata} $\Pi$ is a tuple $(Q, A,\tau,\nu)$ where $Q$ is the set of states  and $A$ alphabet set and ;

The \textit{transition} map $\tau : Q\times A \rightarrow  Q$ and \textit{output} map  $\nu : Q\times A \rightarrow  A$ obey the inductive rules: 
\begin{center}
$q\vert_{\emptyset} = q$    , $q(\emptyset) = \emptyset$  , \\
$q\vert_{x v}= q\vert_{x}\vert_{v}$   ,     $q(x v) = q(x)q\vert_{x}(v)$, 
\end{center}
for $x\in A$ , $v\in A^{*}$ and $q\in Q$ and $q\vert_{x}=\tau(q , x) , q(x)=\nu(q,x)$. 

By the composition low of automatas , define the automata $\Pi^{2}$ over alphabet set $A$ and state set $Q^{2}=Q\times Q$ such that
\begin{center}
$(q_{1}q_{2} , v)\mapsto q_{2}\vert_{q_{1}(v)}q_{1}\vert_{v}$ ,\\    $(q_{1}q_{2} , v)\mapsto q_{1}(q_{2}(v))$
\end{center}
for $v\in A^{*}$ and $q_{1} , q_{2} \in Q$. However by induction one can obtains automaton $\Pi^{*}$ , over alphabet set $A$ with state set $Q^{*}=\bigsqcup_{\mathbb{N}}Q$. The states $Q^{*}$ of the  automata $\Pi^{*}$ with empty set produce  semi-group $G=\langle \Pi^{*} \rangle$.
\end{definition}
Note that the transition $\tau$ and output $\nu$ map can be extented to continuous functions $\hat{\tau}: Q\times A^{\omega}\rightarrow Q$ and $\hat{\nu}:Q\times A^{\omega}\rightarrow A^{\omega}$ repect to the pro-discrete topology.

The automaton $\Pi$ is \textit{invertible} if, for all $q\in \Pi$, the transformation $\tau (q, . ) : A \rightarrow A$ is inertable equivalently $q$ acts by  permutation on $A$. Consequencly the states of an invertable automata $\Pi=\Pi^{-1}$ generate a group is called \textit{automata group} denoted by $G=\langle \Pi \rangle$ which is a subgroup $G\leqslant Aut(A^{*})$ of the atuomorphism group of the rooted tree . An automaton $\Pi$ can be represented by its \textit{Moore diagram}, a directed labeled graph whose vertices are identiﬁed with the states of $\Pi$. For every state $q\in \Pi$ and every letter $x \in A$, the diagram has an arrow from $q$ to $q\vert_{x}$ labeled by $x\vert q(x)$. This graph contains complete information about the automaton, and we will identify the automaton with its Moore diagram.

 The faithful action of a group $G$ on $A^{*}$ (or on $A^{\omega}$) is said to be \textit{self-similar} if for every $g \in G$ and every $x \in A$ there exist $h \in G$ and $y \in A$ such that
\begin{center}
$g(x.w) = y.h(w)$
\end{center}
for every $w \in A^{*}$ (resp. $w \in A^{\omega}$). 
The pair $(h,y)$ is uniquely determined by the pair $(g,x)$, since the action is faithful. Hence we get an automaton with the set of states $G$ and with the output and transition functions 
\begin{center}
$g.x = y.h$
\end{center}
i.e., $y = g(x)$ and $h = g\vert_{x}$. This automaton is called the \textit{complete automaton} of the self-similar action

 In other words, the action of a group $G$ on $A^{*}$ is self-similar if there exists an automaton $(G, A)$ such that the action of $g \in G$ on $A^{*}$ coincides with the action of the state $g$ of the automaton $(G, A)$.\\
The graph of automaton $\mathcal{D}(\Pi) $ is a directed labled graph with vertecs set $V(\mathcal{D}):= Q$ and directed edge $q \rightarrow q\vert_{x}$ labaled $ x \vert q(x)$ equipted with \textit{source} map $s: E(\mathcal{D}) \rightarrow V(\mathcal{D})$ and \textit{target} map $t: E(\mathcal{D}) \rightarrow V(\mathcal{D})$. Let $\Pi $and $\Pi^{\prime}$ be two auromaton on a same alphabet $A$, an \textit{automaton homomorphism} is a a map $f :\Pi \rightarrow \Pi^{\prime}$ such that $f:Q\rightarrow Q$ , $\nu^{\prime}(x , f(q)) = \nu(x , q)$ and $\tau^{\prime}(x , f(q))= f(\tau( x ,f(x)))$ for all $x\in A , q\in Q$. Furthermore , there is labal preserving homomorphism of graphs $\phi:\mathcal{D}(\Pi)\rightarrow \mathcal{D}(\Pi^{\prime})$ if and only if $\phi:\Pi \rightarrow \Pi^{\prime}$ is automaton homomorphism.

\begin{center}
$\psi : G \rightarrow Aut(\mathfrak{M}_{G})\cong Sym(A) \wr G$ , $\psi(g)(m) \in Aut(\mathfrak{M})$ , $m\in\mathfrak{M}$. 
\end{center}

\begin{definition}\label{Nucleus} 
A self-similar action $(G, A)$ is called \textit{contarcting} (or \textit{hyperbolic}) if ther exists a finite set $\mathcal{N}\subset G$ such that for every $g\in G$ there exists $k\in \mathbb{N}$ such that $g\vert_{v}\in \mathcal{N}$ for all words $v\in A^{*}$ of lenght$\geq k$. th minimal set with this property is called the \textit{nucleus} of the sel-similar action.
The nucleus of contracting action is unique and is defiend by:
\begin{center}
$\mathcal{N}:=\bigcup_{g\in G}\bigcap_{n\geq 0}\lbrace g\vert_{v} ; v\in A^{*} , \vert v \vert\geq n\rbrace$
\end{center}
if $g\vert_{v}= g$ for some $v\in A^{*}\setminus\emptyset$, then $g$ belongs to the nuclues by definition.
\end{definition}
\begin{lemma}\cite{nekrashevych2005self}
A sefl-similar action of a group with a generating set $S=S^{-1}$, $1\in S$ is contracting if and only if there exists a finite set $\mathcal{N}$ and a number $k\in \mathbb{N}$ such that 
\begin{center}
$(S\cup \mathcal{N})^{2}\vert_{A^{k}}\subseteq \mathcal{N}$
\end{center}
\end{lemma}
The self-similar action of the group $G$ on a rooted tree  $A^{*}$ is spherically transitive and faithful if the action of $G$  on each level $A^{n}$ is transitive and faithful. Moreover, the action extend to a continuous transitive faithful action on the boundary $\partial\mathcal{T}$ (see \cite{bass2006cyclic}).
\begin{definition}\cite{bartholdi2003fractal}\cite{nekrashevych2005self}
A self-similar action is said to be \textit{self-replacing} (\textit{recurrent} or \textit{fractal}) if it is transitive on the ﬁrst level $A$ of the tree $ A^{*}$ and the map $g \mapsto g\vert_{x}$ from the stabilizer $Stab_{G}(x)$ to the group $G$ is surjective for some
(every) letter $x \in A$. It can be shown that a self-replicating group acts transitively on $X_{n}$
for every $n \geq 1$. It is also easy to see (\cite{nekrashevych2005self}, Proposition 2.11.3) that if a finitely generated
contracting group is self-replicating then its nucleus$\mathcal{N}$ is a generating set.
\end{definition}
It is easy to see that when the action is recurrent the transivity on the first level $A$ results the transivity on all levels $A^{n}$ of rooted tree $A^{*}$. 
\begin{remark}\cite{nekrashevych2005self}
 Even though it is usually challenging to determine whether the automata
generate a contracting group. It takes polynomial time to determine the nucleus of self-similar group action.
\end{remark}

\subsubsection{Schreier graphs}
\begin{definition}\label{definition Schreier graph}
Let $G$ be a group generated by a ﬁnite set $S$ and let $H$ be a subgroup of $G$. The \textit{Schreier graph} $\Gamma(G,S,H)$ of the group $G$ is the graph whose vertices are the right cosets $G/H = \lbrace Hg : g \in G\rbrace$, and two vertices $Hg_{1}$ and $Hg_{2}$ are adjacent if there exists $s \in S$ such that $g_{2} = g_{1}.s$ or $g_{1} = g_{2}.s$. 
suppose the finitely generated group $G$ acts on a set $M$. Then the corresponding \textit{Schreier graph}  $\Gamma(G,S,M)$ is the graph with the set of vertices $M$ and set of arrows $S \times M$, where the arrow $(s,v), v\rightarrow s(v)$ starts in $v$ and ends in $s(v)$. The \textit{simplicial Schreier graph} $\overline{\Gamma}(G,S,M)$ remembers only the vertex adjacency:(see figur) its set of vertices is $M$ and two vertices are adjacent if and only if one is an image of another under the action of a generator $s \in S$. If $(G,A)$ is a self-similar action and $G$ is generated by a ﬁnite set $S$, then we get a sequence $\Gamma_{n} = \Gamma(S,A^{n})$ of Schreier graphs of the action of G on the levels $A^{n}$ of the rooted tree $A^{*}$ (and the sequence $\overline{\Gamma_{n} }= \overline{\Gamma}(S,A^{n})$ of the respective simplicial Schreier graphs). If $(G,A)$ is generated by a ﬁnite automaton $\Pi (Q,A)$, then the graphs$ \Gamma(Q,A^{n})$ coincide with the dual Moore diagrams of the automata $\Pi^{\prime} (Q,A)$ .

If the action of $G$ is transitive, then the Schreier graph $\Gamma(G,S,M)$ is isomorphic to the Schreier graph $\Gamma(G,S,St_{G}(x))$ of the group where  $St_{G}(x)$ is the stablizer of some point $x \in A$.

Let $G$ be a self-similar group generated by a finite set $S$. The sets $A^{n}$ are invariant under the action of $G$, and we denote the associated Schreier graphs by $\Gamma_{n} = \Gamma_{n}(G,S, A^{n})$. For a point $\xi \in A^{\omega}$ we consider the action of the group $G$ on the $G$-orbit of $\xi$, and the associated Schreier graph is called \textit{orbital Schreier graph} denoted $\Gamma_{\xi} = \Gamma_{\xi}(G,S)$. For every $\xi \in A^{\omega}$ we have $St_{G}(\xi) = \bigcap_{n\geq1} St_{G}(\xi_{ n})$, where $\xi_{n}$ denotes the preﬁx of length $n$ of the inﬁnite word $\xi$. The connected component of the rooted graph $(\Gamma_{n} ,\xi_{n})$ around the root $\xi_{n}$ is exactly the Schreier graph of $G$ with respect to the stabilizer of $\xi_{n}$. It follows immediately that the graphs $(\Gamma_{n} ,\xi_{n})$ converge to the graph $(\partial\Gamma ,\xi)$ in the pointed Gromov-Hausdorﬀ topology. 

\end{definition}
\begin{figure}[h!]
  \begin{subfigure}[b]{0.45\linewidth}
  \centering
    \includegraphics[width=0.9\textwidth]{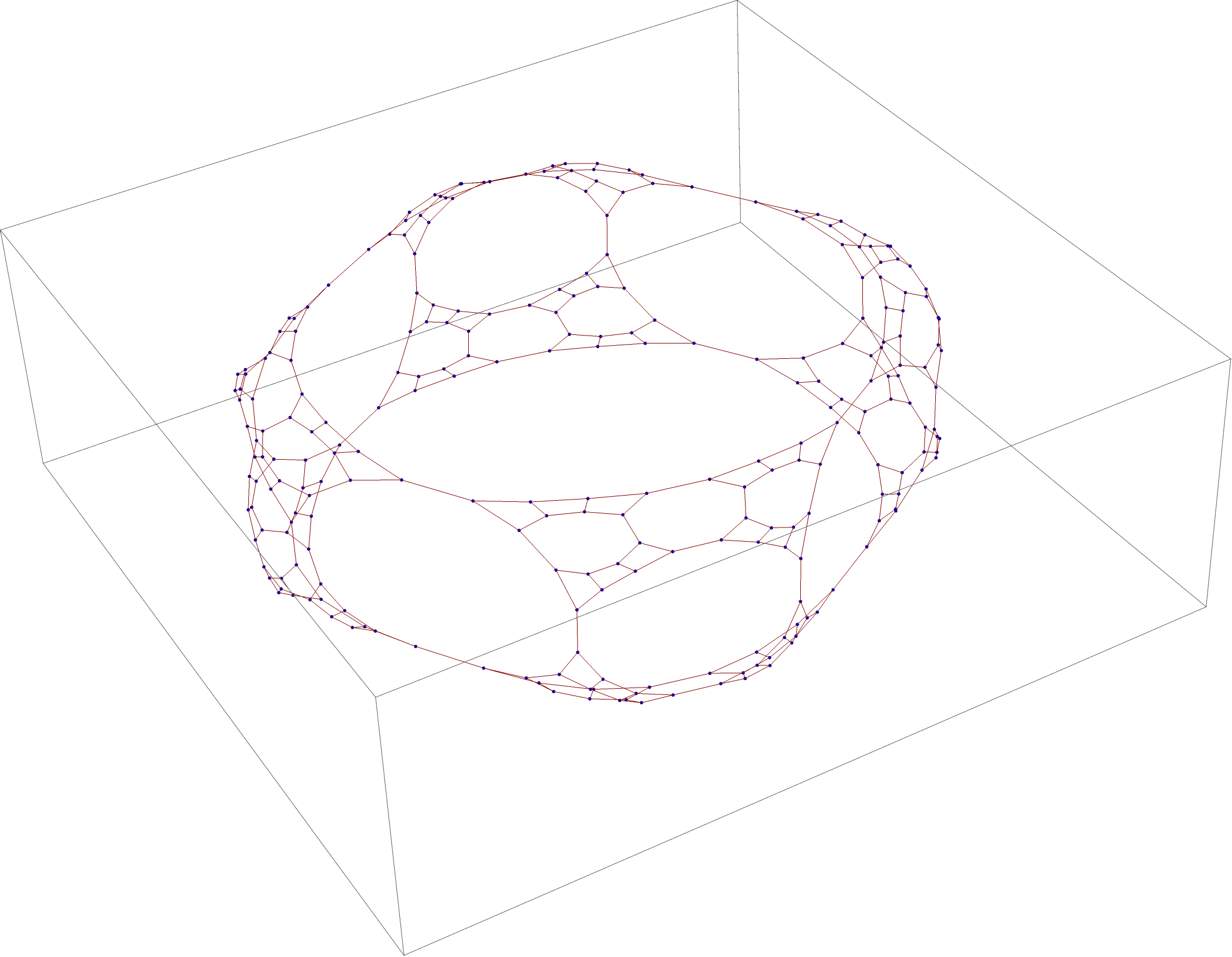}
    \caption{
 The Schreier graph level 8 }
  \end{subfigure}
\begin{subfigure}[b]{0.45\linewidth}
\centering
    \includegraphics[width=0.9\textwidth]{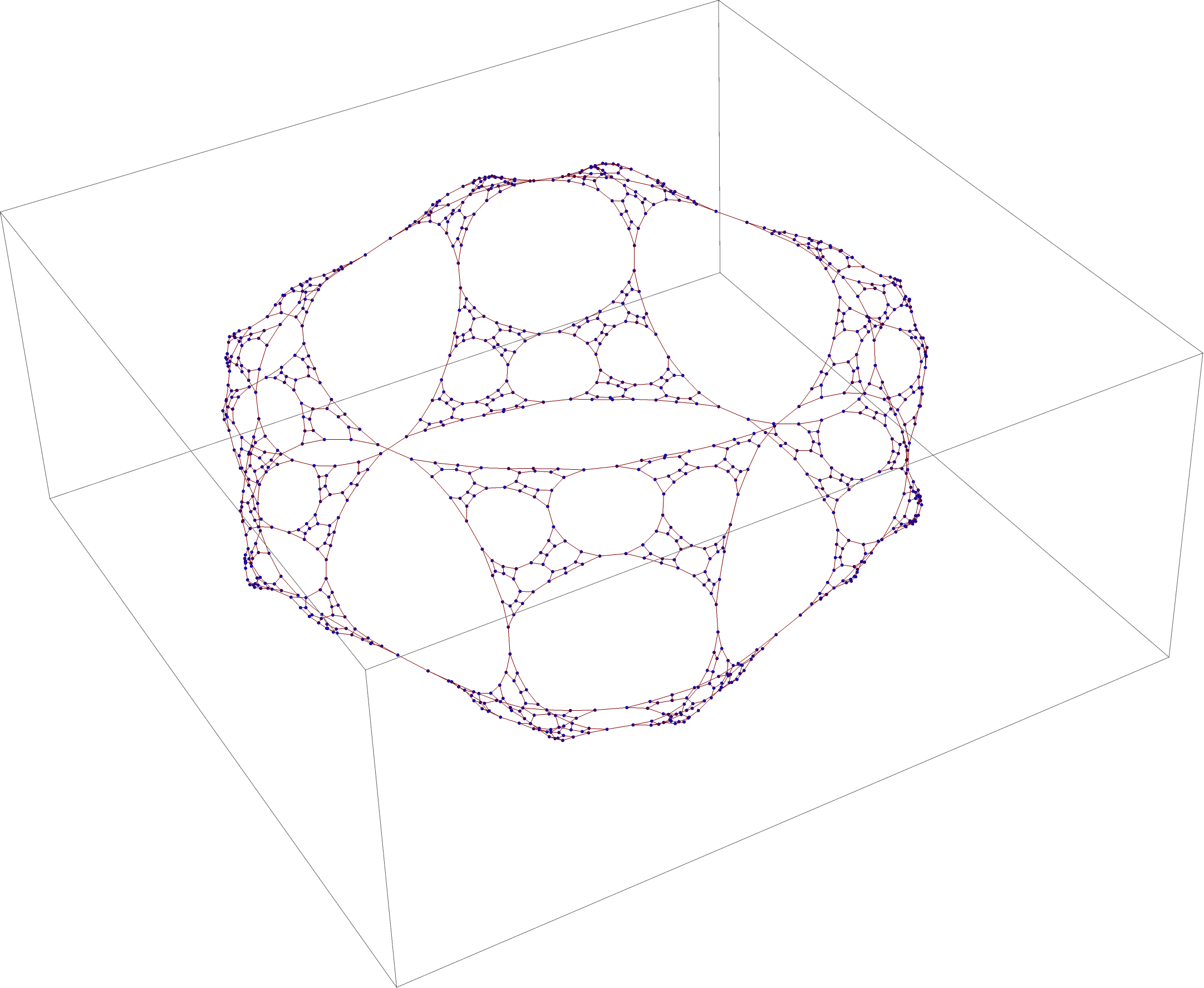}
     \caption{
 The Schreier graph level 10  }
  \end{subfigure}
  
\begin{subfigure}[b]{0.45\linewidth}
\centering
    \includegraphics[width=1\textwidth]{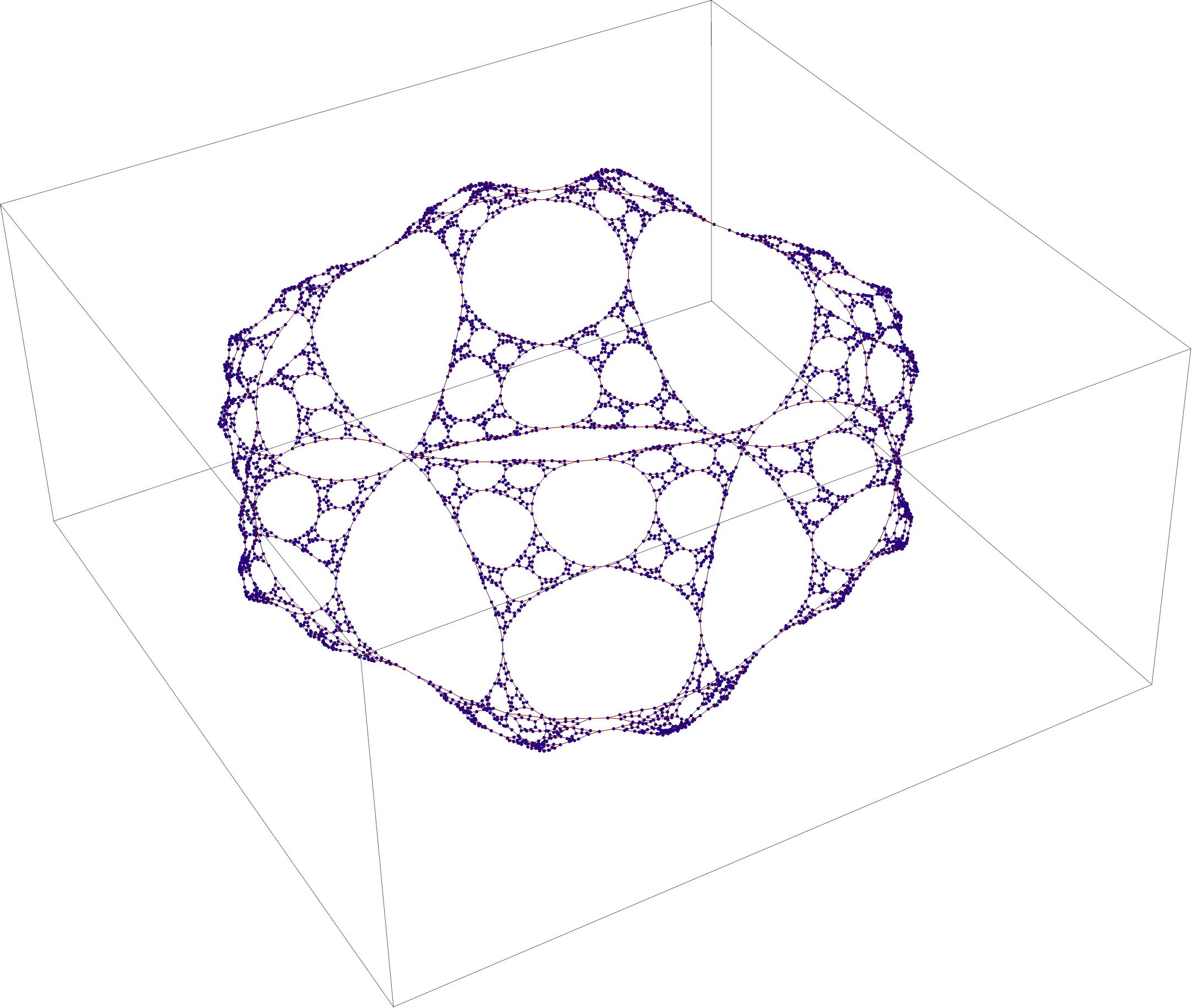}
     \caption{
 The Schreier graph level 12  }
  \end{subfigure}
  \caption
{The Schreier graph of $IMG(\dfrac{z-1}{z+1}^{2})\cong$ group generated with the wreath recursion $a=\sigma(c,c) , b=\sigma (b,a) ,c=(c,c)$.}

  \label{2853 figure}
\end{figure}
\subsubsection{Limit space and the ends of trees}
The limit space space can be defined in different ways: as a quotient of the Cantor set by an
asymptotic equivalence relation (Definition 9.1), as a inverse limit of finite Schreier graphs (Theorem 9.6)
or as the boundary of a $\delta$-hyperbolic graph.

Assume  $G$ to be an infinite finitely generated group with finite generating
set $S$. There is a natural geometric left action of $G$ on the Cayley graph
$\Gamma(G,S)$ equipped with word metric. One can obtain a natural compactification $\overline{\Gamma}(G,S)$ of $X$ by corresponding to $X$ the inverse limit of a system of closed balls in the Cayley graph concentrated on an arbitrary point $x_{0}\in \Gamma(G,S)$ with radius $r$ and projection map  $p_{r}: \overline{B}_{\Gamma}(x_{0},r^{\prime})\rightarrow \overline{B}_{\Gamma}(x_{0},r)$ for $r^{\prime}\geq r$. There is natural continuous (respect to the inverse topology) bijection $\phi(x_{0}):\overline{\Gamma}(G,S)\rightarrow \varprojlim\overline{B}_{\Gamma}(x_{0},r)$ which induces a topology on $\overline{\Gamma}(G,S)$ sometimes called \textit{cone} topology \cite{bridson2013metric}. The induced topology coincides with the topology of uniform convergence on compact subsets. Thus $\varprojlim\overline{B}_{\Gamma}(x_{0},r)$ is a compactification of the Cayley graph $\Gamma$; it may be viewed as a subspace of the set of continuous maps, which associated to each point $\xi \in\partial \Gamma=\overline{\Gamma}\setminus \Gamma$, a geodesic ray $c:[0,\infty) \rightarrow \Gamma(G,S)$ that issues from $x_{0}$ that converge to  is referred to \textit{ideal point} or \textit{point at infinity}. That $\partial \Gamma$ is a closed subspace of $\overline{\Gamma}$ and $\overline{\Gamma}$ is compact. Notice that the asymptotic geometry of Cayley graph for group $G$ is independent of the choice of the generator se, so the notation $\partial \Gamma(G,S)$ can relax to $\partial G$. 

Let $G=\langle \Pi \rangle$ be a automatic group generated by the states of an automata over alphabet set $A$, which acts self-similary on spherically homogeneous rooted tree $A^{*}$. The \textit{ends}(boundary at infinity) of the tree  denoted by $A^{-\omega}$ consists of the left infinite sequence of symbols. The group action admits the asymptotic equivalence relation $\sim$ on $A^{-\omega}$. Two left-infinite sequences $\cdots x_{2}x_{1}$ and $\cdots y_{2}y_{1}$ are equvalence if there is left infinite path in the Moor diagram labled $\cdots (x_{2},y_{2})(x_{1},y_{1})$. The limit space of the self-similar action denoted by $\mathcal{J}_{G}$ is the quotient of the
topological space $A^{\omega}$ by the asymptotic equivalence relation $\sim$.

H.Bass \cite{bass2006cyclic} shows that if an element $g\in Aut(A^{*})$ acts spherically transitively and faithful on each level of the tree, this action extends to a continuous faithful transitive on the boundary of the tree.

\begin{theorem}\cite{nekrashevych2005self}
 Let $G$ be a finitely generated group with a recurrent contracting action; then the limit $G$-space $\mathcal{J}_{G}$ is connected and locally connected. 
\end{theorem}
One can immediately reach If the action of a finitely generated group is contracting and recurrent; then the limit space is path-connected and locally path-connected.

Despite the homeomorphic action of the group $G$ on the boundary
$ \partial  \Sigma(G, S,X)$, the limit space $\mathcal{J}_{G}$ is a hyperbolic boundary. For a definition of hyperbolic boundary and more details, the reader is referred to \cite{gromov1999topological}\cite{kapovich2002boundaries}\cite{bridson2013metric}.

\begin{definition}
Let $G$ is a groups with generating set $S$ with self-similar action on $A$.
The \textit{self-similarity graph} $\sum(G,S,A)$ is a graph with $A^{*}$ as vertices set and to vertices $v_{1},v{2}\in A^{*}$  are adjacent  if they are connectet either by a vetrical edge if $v_{1}=x v_{2}$ for some $x\in A$, or horizontal edge if $s(v_{i})=v_{j}$ , $i,j\in {1,2}$ for some $s\in S$. 
\end{definition}
 The following theorem determine some important aspect of the limit space when the group $G$ acts by contraction on the rooted tree $A^{*}$. 
\begin{theorem}\cite{nekrashevych2005self}
The limit space $\mathcal{J}_{G}$ of a contracting action of a finitely generated
group $G$ is homeomorphic to the hyperbolic boundary $\partial  \Sigma(G, S,X)$ of the self-similarity
graph $\Sigma(G, S,X)$. Moreover, there exists a homeomorphism $F : \mathcal{J}_{G}\rightarrow
\partial  \Sigma(G, S,X)$, such that $D = F \circ \pi$, were $\pi : X^{-\omega}\rightarrow \mathcal{J}_{G}$ is the canonical projection
and $D :X^{-\omega}\rightarrow\partial  \Sigma(G, S,X)$ carries every sequence $\cdots x_{2}x_{1} \in X^{-\omega}$ to its limit
\begin{center}
$\underrightarrow{lim}_{\infty}x_{n}\cdots x_{2}x_{1}\in \partial  \Sigma(G, S,X) $
\end{center}
\end{theorem}

The asymptotic equivalence relation is closed and has finite
equivalence classes. The limit space $\mathcal{J}_{G}$ is a compact metrizable with finite topological dimension $\leq\vert \mathcal{N}\vert$.
 \begin{example}[Sierpinski gasket]
 Let $A=\lbrace0,1,2\rbrace$ be the alphabet set and consider the group $G$ generated by wreath recurssion:
 \begin{center}
$a=\sigma_{1}(e, a, e) , b=\sigma_{2}(e , e , b) , c=\sigma_{2}(c , e , e) , e=(e,e,e)$\\
 $\sigma_{1}=(0,2) , \sigma_{2}=(0,1) , \sigma_{3}=(1,2)$
\end{center}
where $\sigma_{i} , i=1,2,3$ are permutation of $A$.
\begin{figure}[h]
 \begin{subfigure}[b]{0.45\linewidth}
 \centering
    \includegraphics[width=0.7\textwidth]{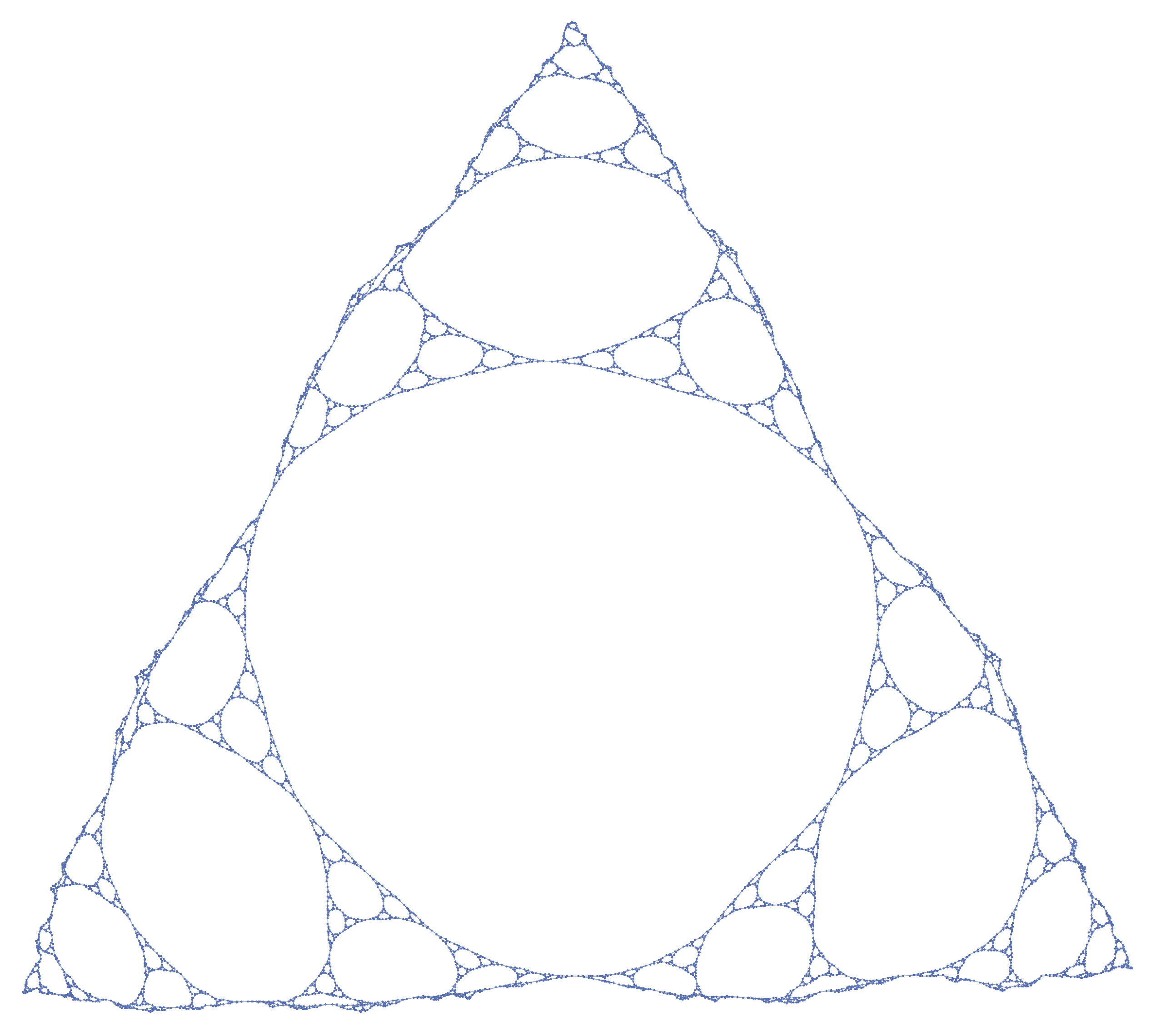}
    \caption{The tail graph level 8}
  \end{subfigure}
\begin{subfigure}[b]{0.45\linewidth}
\centering
    \includegraphics[width=0.8\textwidth]{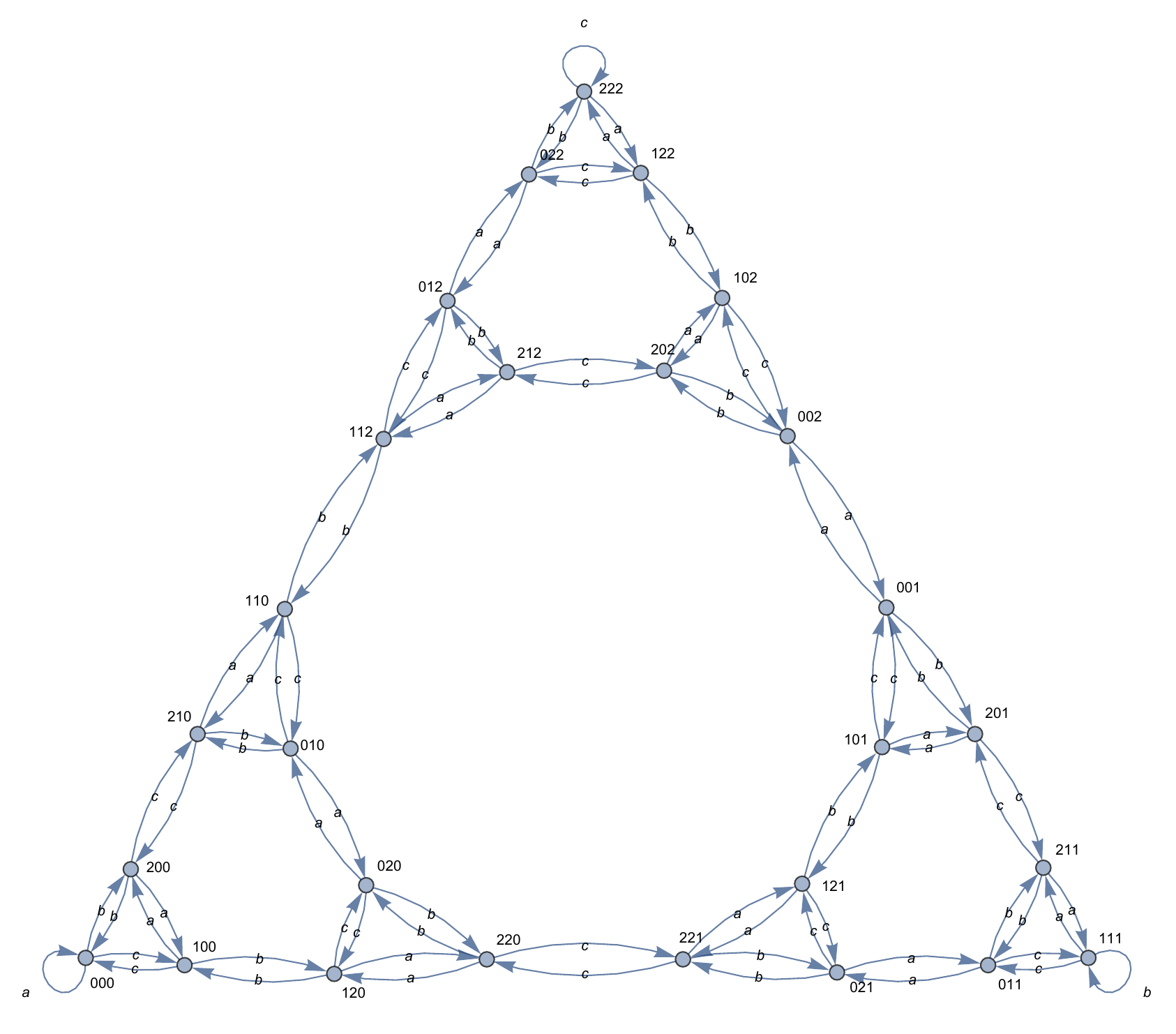}
     \caption{The Schreier graph level 3}
  \end{subfigure}
  \caption{  Sierpinski group action on 3-regular rooted tree at level 8-th }
  \label{Sierpinski group}
   \end{figure}
\end{example}

\section{Approximating with the Schreier graphs}\label{result}
In this section, we use the proposed program to explore the incidence structure of Schreier graphs related to groups generated with finite automaton and assemble the approximation of the limit space $\mathcal{J}_{\langle \Pi \rangle}$ by their Schreier graphs. The program is applied to a set of automata which contains the first Grigorchuk group, Basilica group, the free group $\mathbb{
F}_{3}$ of rank 3, Sierpinski group, the lamplighter group, the mother groups, the longe range group, Hanoi tower group, and some automatic groups that are represented in \cite{grigorchuk2014self}\cite{bondarenko2008classification}.
\begin{itemize}

\item[1- ]The automata by the wreath recurssion $a = \sigma(b, c)
b = \sigma(c, b)
c = (a, a)$get the free group of rank tree first time introduced by Aleshin
The associated limit space has not have a finite Hausdorff dimension (Figure \ref{free group}) . 
\item[2- ]
Let $G=\langle \Pi \rangle$ be an automata group generated by Wreath recursion:
\begin{center}
$  a = \sigma(b, c), b =\sigma(c, b), c =(a, a) $.
\end{center}
Aleshin \cite{aleshin1972finite} introduced this automaton. \cite{nekrashevych2010free} has proved that group $G$ is a free group of
rank three. The automaton is the smallest automaton among all automata
over a 2-letter alphabet, generating a free non-abelian group.

\begin{figure}[h]
\centering
    \includegraphics[width=0.8\textwidth]{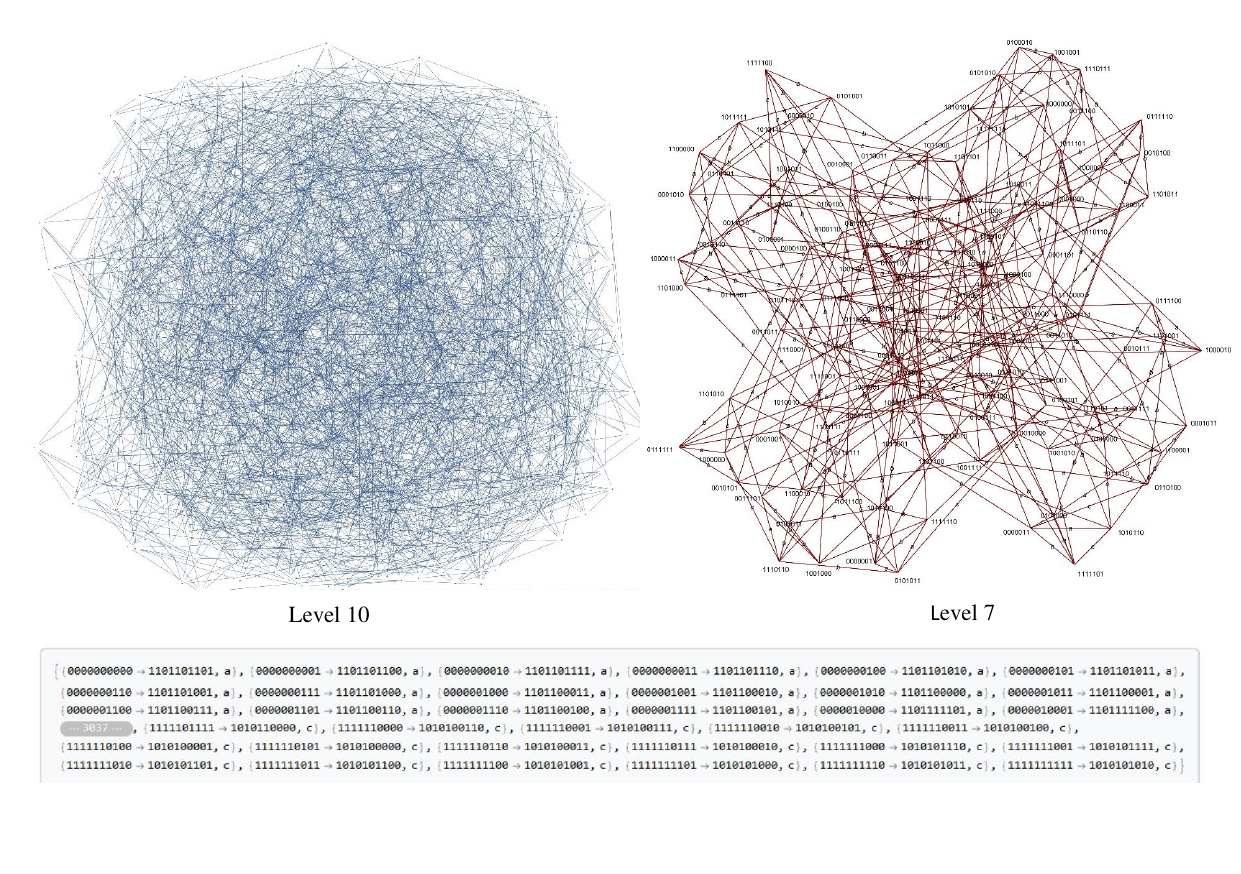}
     \caption{The Schreier graph of the $\mathbb{F}_{3}$ free group of rank three at level $n=7,10$}
     \label{free group}
\end{figure}

\item[3- ]
The automata 882, produced by the wreath recursion \cite{bondarenko2008classification}:

\begin{center}
$ a = \sigma(c, c), b = (b, c), c = (b, a)$.
\end{center}
The element $ (ca^{-1}cb^{-1})^{2}$  stabilizes the vertex $00$ and its section at
this vertex is equal to $ca^{-1}cb^{-1}$. Hence, $ca^{-1}cb^{-1}$ has infinite order.
\begin{figure}[h]
\centering
\begin{subfigure}[b]{0.7\linewidth}
\centering
    \includegraphics[width=1\textwidth]{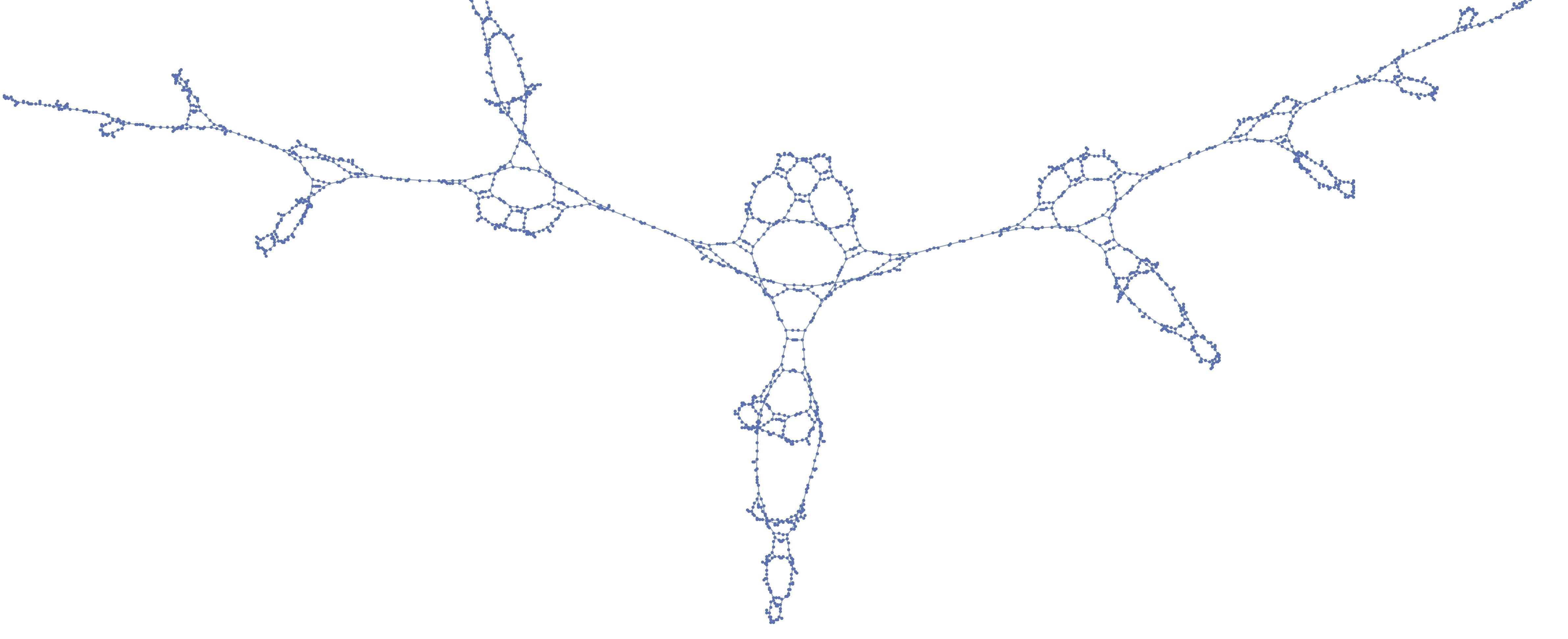}
     \caption{The Schreier graph}
  \end{subfigure}
  \begin{subfigure}[b]{0.8\linewidth}
\centering
    \includegraphics[width=1\textwidth]{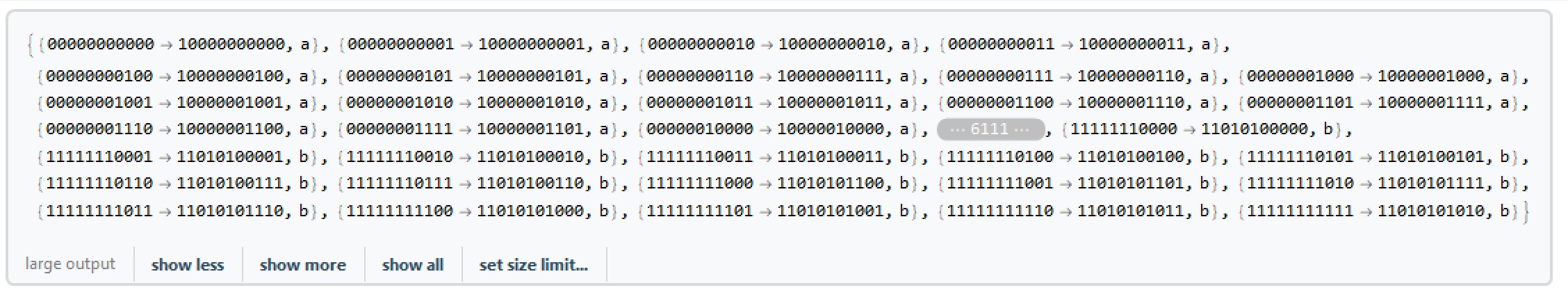}
     \caption{A small part of adjacency computation}
  \end{subfigure}
  \caption{The action of automata generated by wreath recurssion $ a = \sigma(c, c), b = (b, c), c = (b, a)$  on level 11-th binary rooted tree }
  \label{882 figure}
\end{figure}

\item[4- ]
the case of self-similar actions of the free abelian group
$\mathbb{Z}^{n}$ the limit space $\mathcal{X}_{\mathbb{Z}^{n}}$ is homeomorphic to $\mathbb{R}^{n}$  and the tile $\mathcal{T}$  is an integral self-affine tile (Figure \ref{Z2}),
   \begin{figure}[h]
  \begin{subfigure}[b]{0.4\linewidth}
  \centering
    \includegraphics[width=0.8\textwidth]{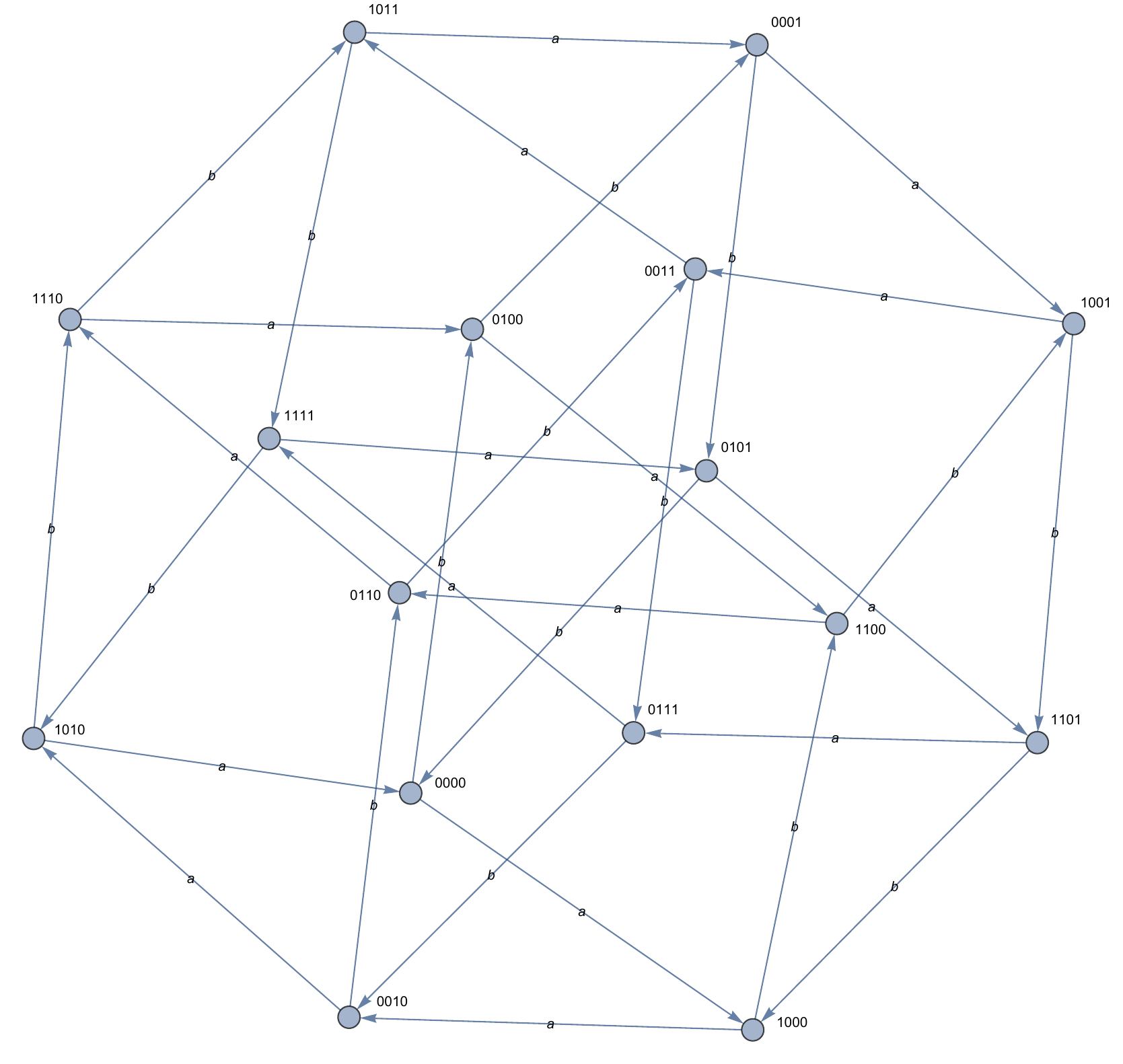}
    \caption{
 The Schreier graph level 4 }
  \end{subfigure}
\begin{subfigure}[b]{0.5\linewidth}\centering
    \includegraphics[width=0.8\textwidth]{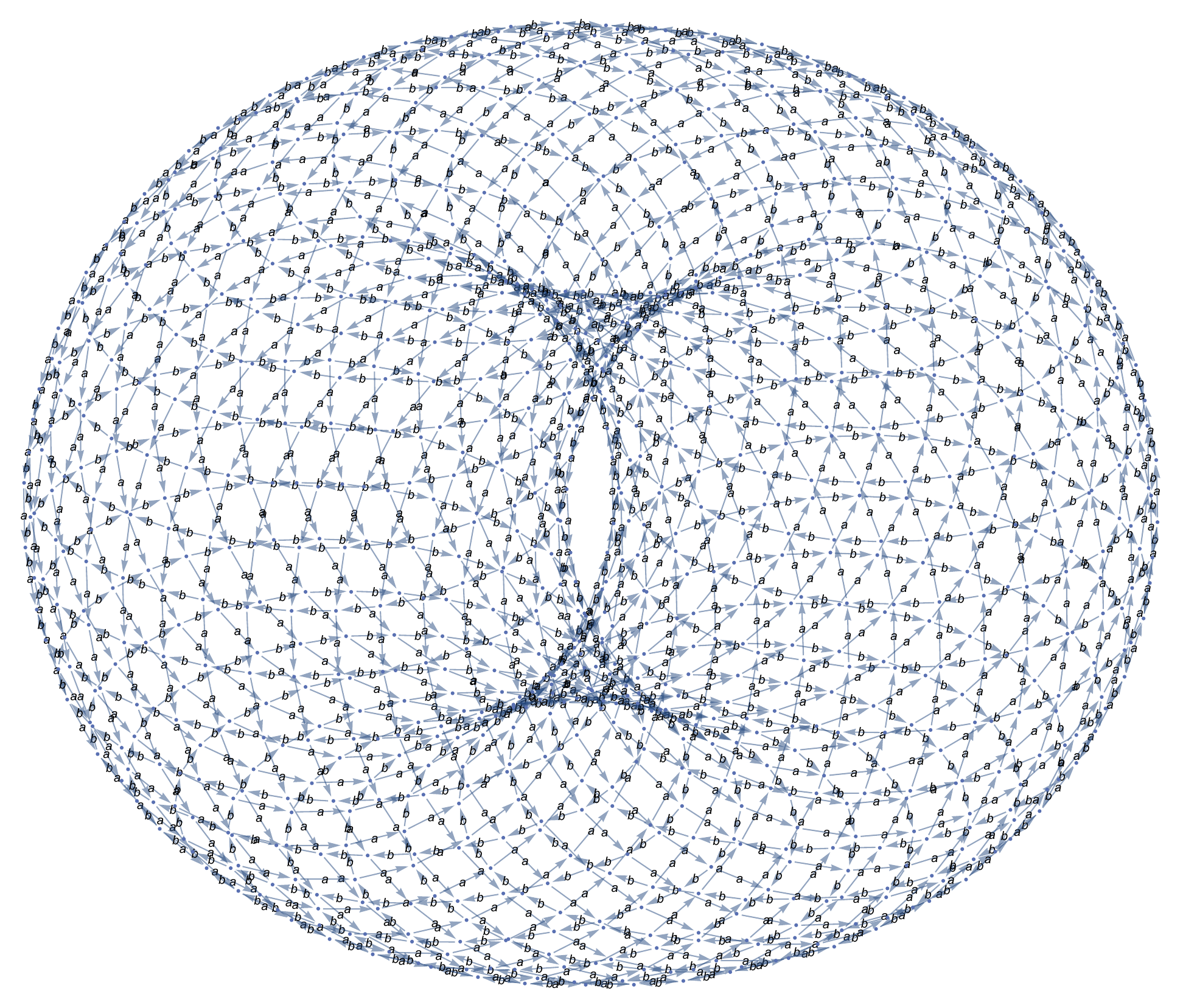}
     \caption{
 The Schreier graph level 10  }
  \end{subfigure}
  \caption
{The self-similar contracting action of $\mathbb{Z}^{2}$ on  binary rooted tree, generated with wreath recurssion $a=\sigma (e, b),b=(a,a),e=(e,e)$ }
  \label{Z2}
\end{figure}

\item[5- ]
Consider the wreath recursion:
 
\begin{center}
$  a = \sigma (b, b), b = (b, c),
c = (b, a)$
\end{center}
Let $x = bc $ and $y = ca$. Since all generators have order $2$, the index
of the subgroup $H = \langle x, y\rangle$ in $G_{878}$ is $2$, $H$ is normal and $G_{878}\cong C_{2} \ltimes H$, where $C_{2}$ is generated by $c$. The action of $C_{2} $ on $H$ is given by
$x^{c} = x^{-1}$ and $y^{c} = y^{-1}$. We have $x = bc = (1, ca) = (1, y)$ and $y =
ca = \sigma (ab, 1) = \sigma (y^{-1}x^{-1}, 1)$. An isomorphic copy of $H$ is obtained by
exchanging the letters $0$ and $1$, yielding the wreath recursion (Figure \ref{878 figure})

\begin{center}
$x = (y, 1)
$ and $y = \sigma (1, y^{-1}x^{-1})$. 
\end{center}
The last recursion defines $IMG(1 - \dfrac{1}{z^{2}}
 )$ for more details the reader are refered to\cite{bartholdi2006thurston}. Thus, $\langle\Pi_{878}\rangle\cong C_{2} \rtimes IMG(1 -\dfrac{1}{z^{2}}
 )$.
%   \begin{figure}[h]
%  
%  \begin{subfigure}[b]{0.45\linewidth}
%  \centering
%    \includegraphics[width=1.1\textwidth]{automata878level4.pdf}
%    \caption{level 4 }
%  \end{subfigure}
%\begin{subfigure}[b]{0.45\linewidth}
%\centering
%    \includegraphics[width=1.4\textwidth]{878level11.pdf}
%     \caption{level 11}
%  \end{subfigure}
%  \caption{The Schreier graph of $IMG(1 -\dfrac{1}{z^{2}})$ generated by wreath recurssion $  a = \sigma (b, b), b = (b, c),
%  c = (b, a)$}
%\end{figure}
\begin{figure}[h!]
  \begin{subfigure}[b]{0.45\linewidth}
  \centering
 \includegraphics[width=1.1\textwidth]{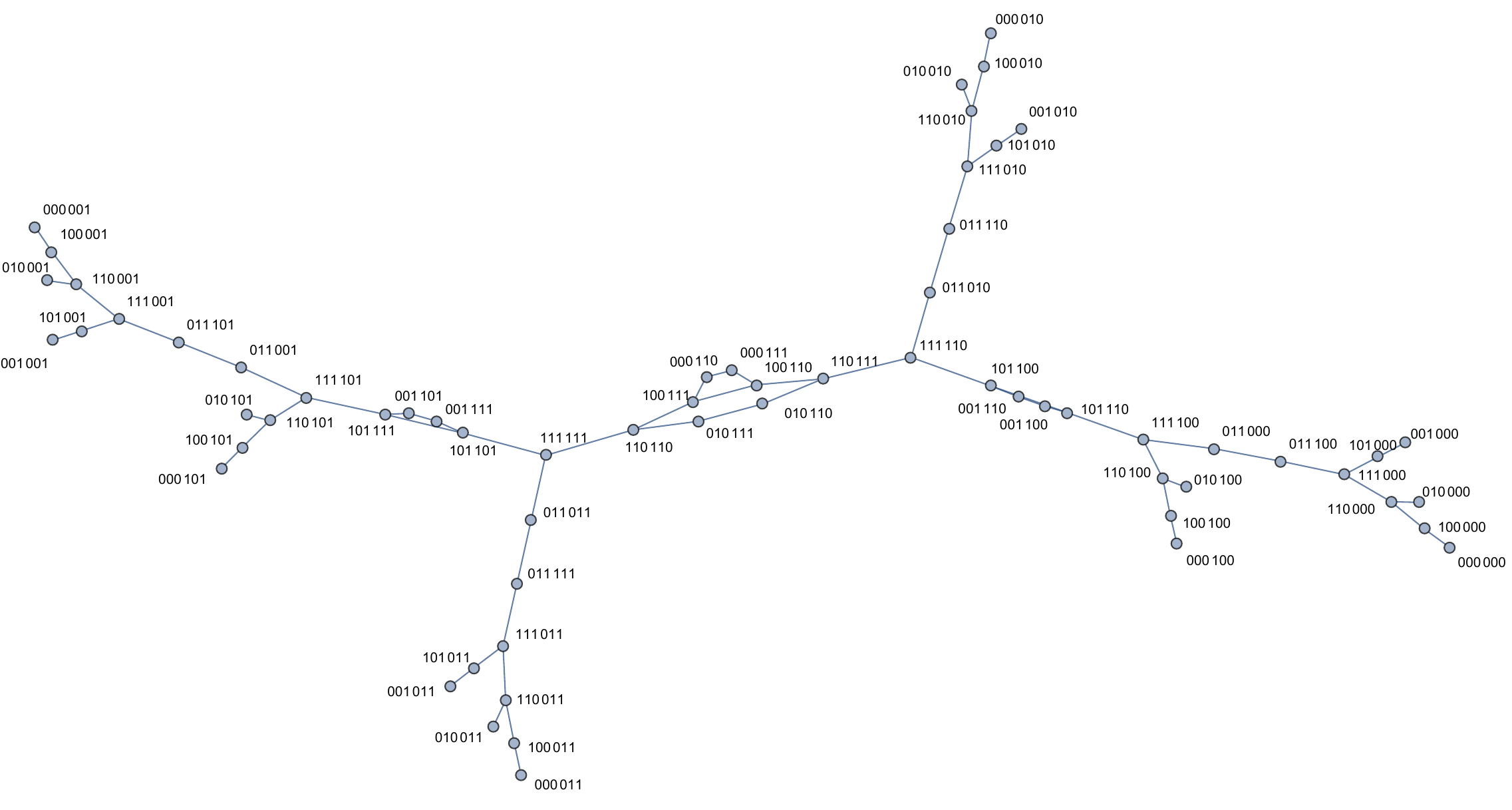}
    \caption{ level 6 }
  \end{subfigure}
\begin{subfigure}[b]{0.45\linewidth}
\centering
    \includegraphics[width=1.1\textwidth]{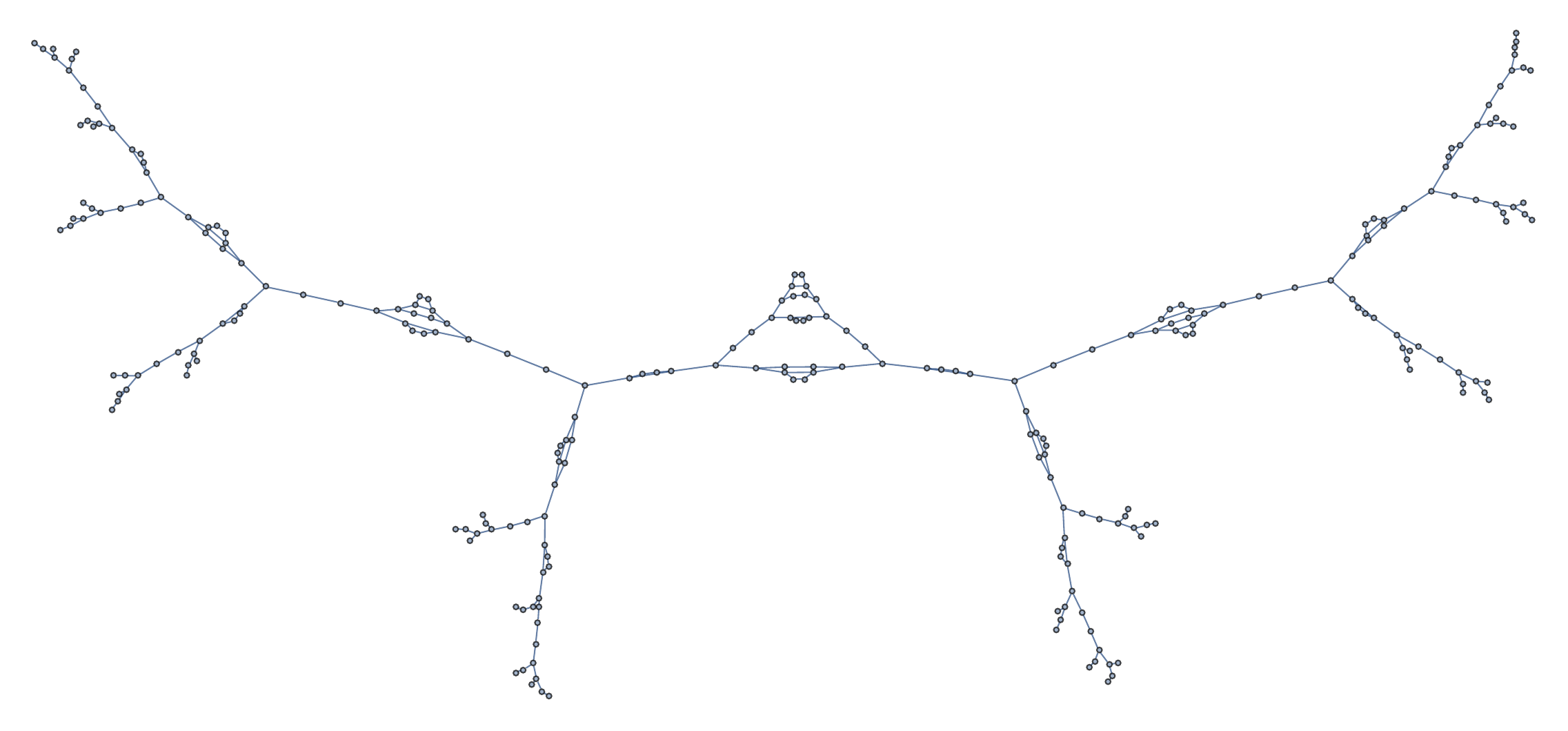}
     \caption{ level 8  }
  \end{subfigure}
  
\begin{subfigure}[b]{0.45\linewidth}
\centering
    \includegraphics[width=1\textwidth]{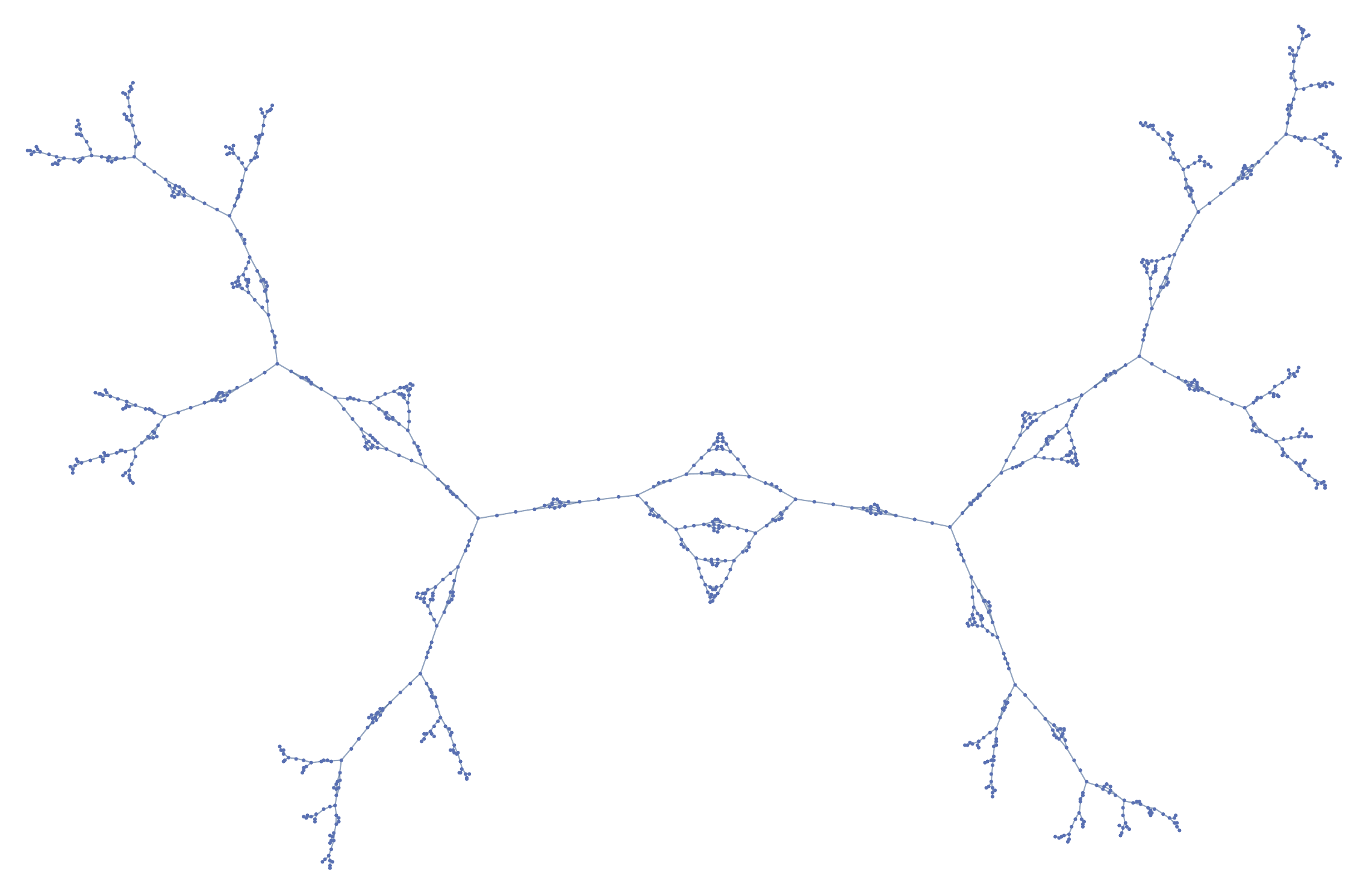}
     \caption{ level 10}
  \end{subfigure}
  \begin{subfigure}[b]{0.45\linewidth}
\centering
    \includegraphics[width=1\textwidth]{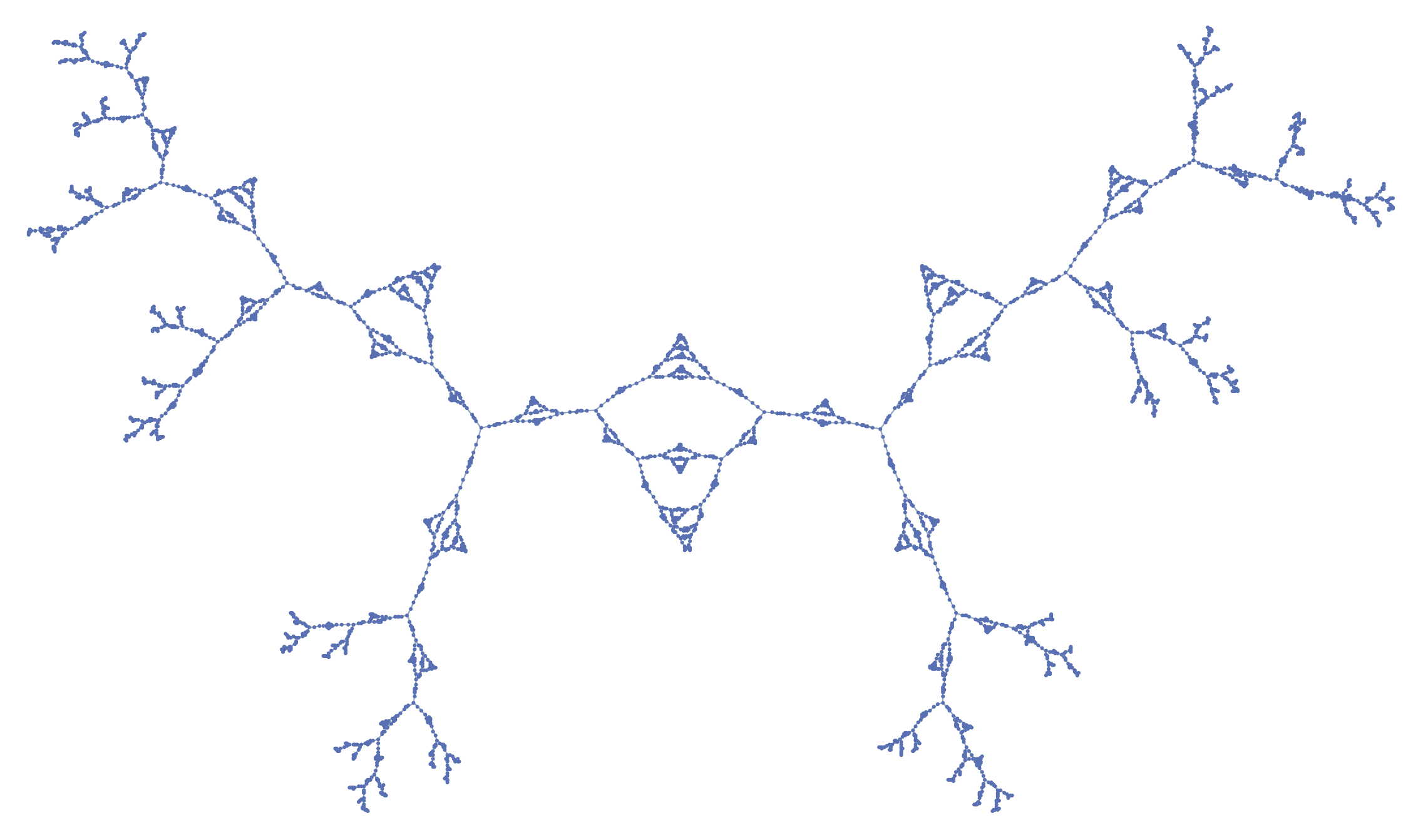}
     \caption{ level 12}
  \end{subfigure}
  \caption
{The simplicial Schreier graph of $IMG(1 -\dfrac{1}{z^{2}})$ generated by wreath recurssion $  a = \sigma (b, b), b = (b, c),
c = (b, a)$.}
  \label{878 figure}
  \end{figure}

\item[6-]\label{half-Basilica group}
Let $G$ be a group generated by a finite state automata with wreath recursion:
\begin{center}
$a=\sigma(b,b) , b=(c,b) , c=(c,a)$
\end{center}
The action of group $G$ on the binary tree is conjugated with the half-Basilica group $C_{2}\ltimes IMG(z^{2}-1)$ (see figure \ref{half-Basilica group figure}).
 \begin{figure}[h!]
\centering
    \includegraphics[width=0.6\textwidth]{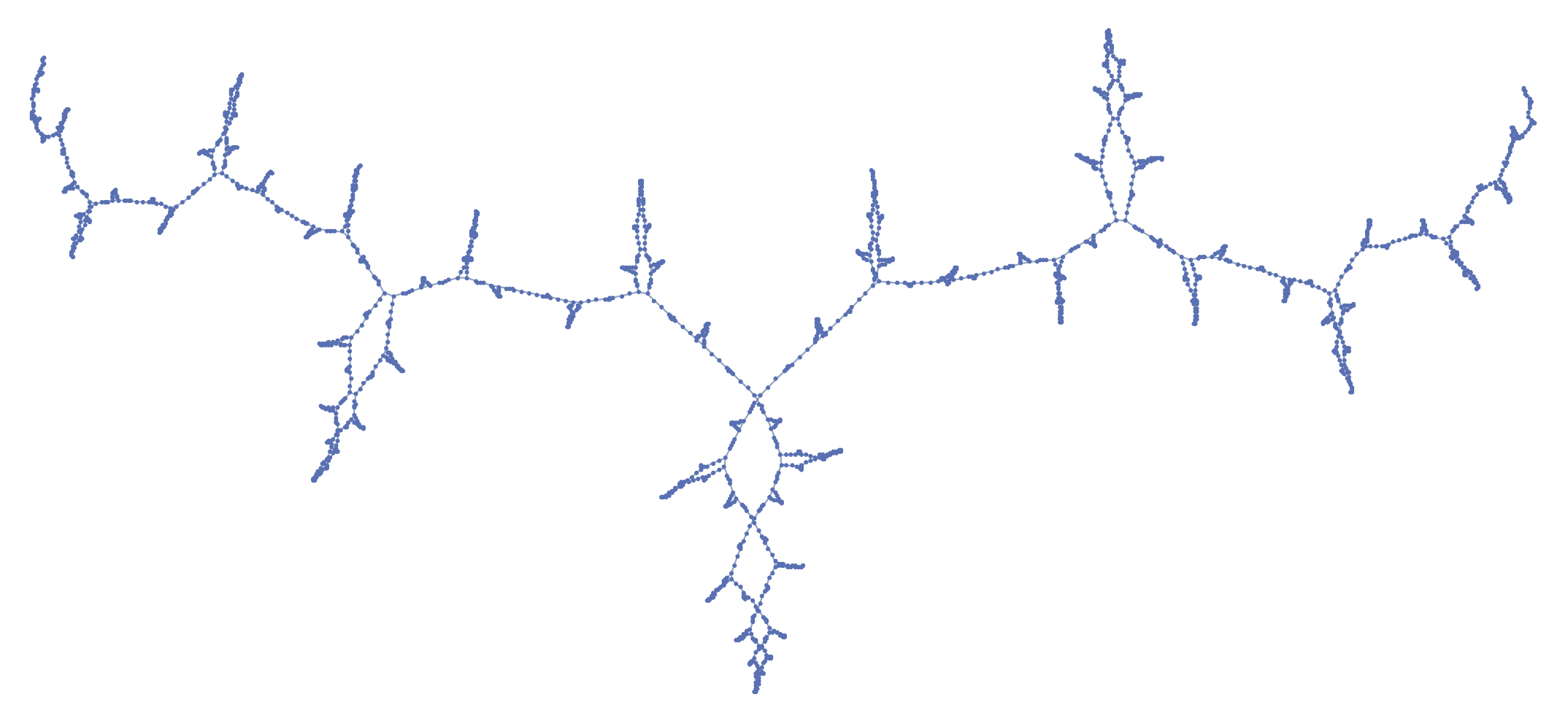}
  \caption{The Schreier graph generated by automata $a=\sigma(b,b) , b=(c,b) , c=(c,a)$ isomorphic to the half-Basilica group $C_{2}\ltimes IMG(z^{2}-1)$}
  \label{half-Basilica group figure}
\end{figure}
\item[7-]\cite{amir2013amenability}
The long-range group, an interesting group generate by the wreath recursion
$a = (a ,b), b=\sigma(b,1) $. where $\sigma=(01)$ is permutation of $\lbrace 0,1\rbrace$
\begin{figure}[h]
\centering
    \includegraphics[width=0.5\textwidth]{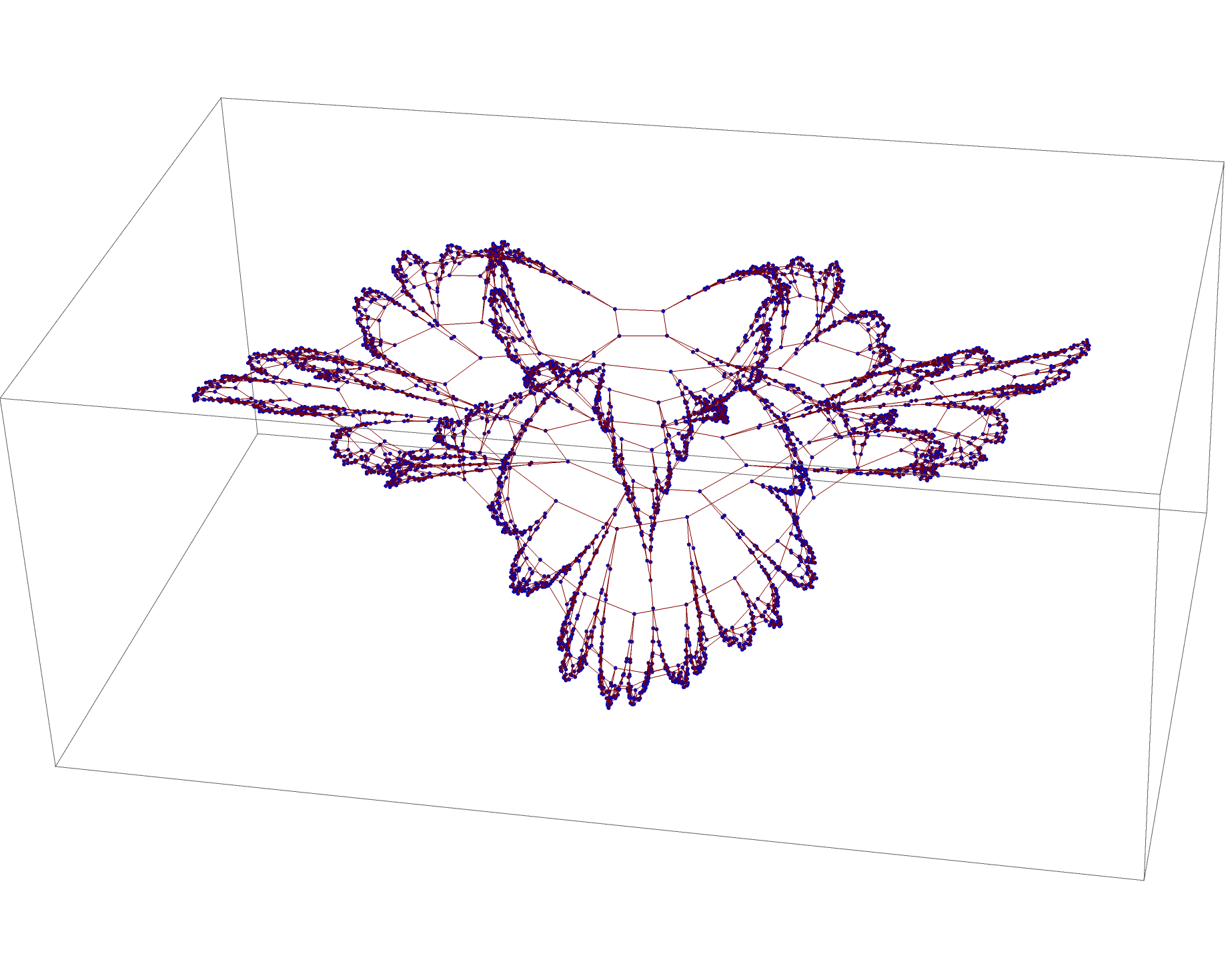}
     \caption{The Schreier graph of the long-range group at level 12 }
     \label{long-range}
\end{figure}
The Schreier graphs of long-range group (Figure \ref{long-range}) were studied by
Benjamini and Hoffman (2005)\cite{benjamini2003omega} in the context of long-range percolation theory.

\item[8-] The \textit{Hanoi tower}s group, the group of possible moves in the Hanoi Towers game on three
pegs, a game introduced by ` Edouard Lucas in 1883. Its level Schreier graphs are discrete
Sierpinski gaskets (see\cite{grigorchuk2007schreier},\cite{grigorchuk2006asymptotic}).
\begin{figure}[h]
\centering
    \includegraphics[width=0.45\textwidth]{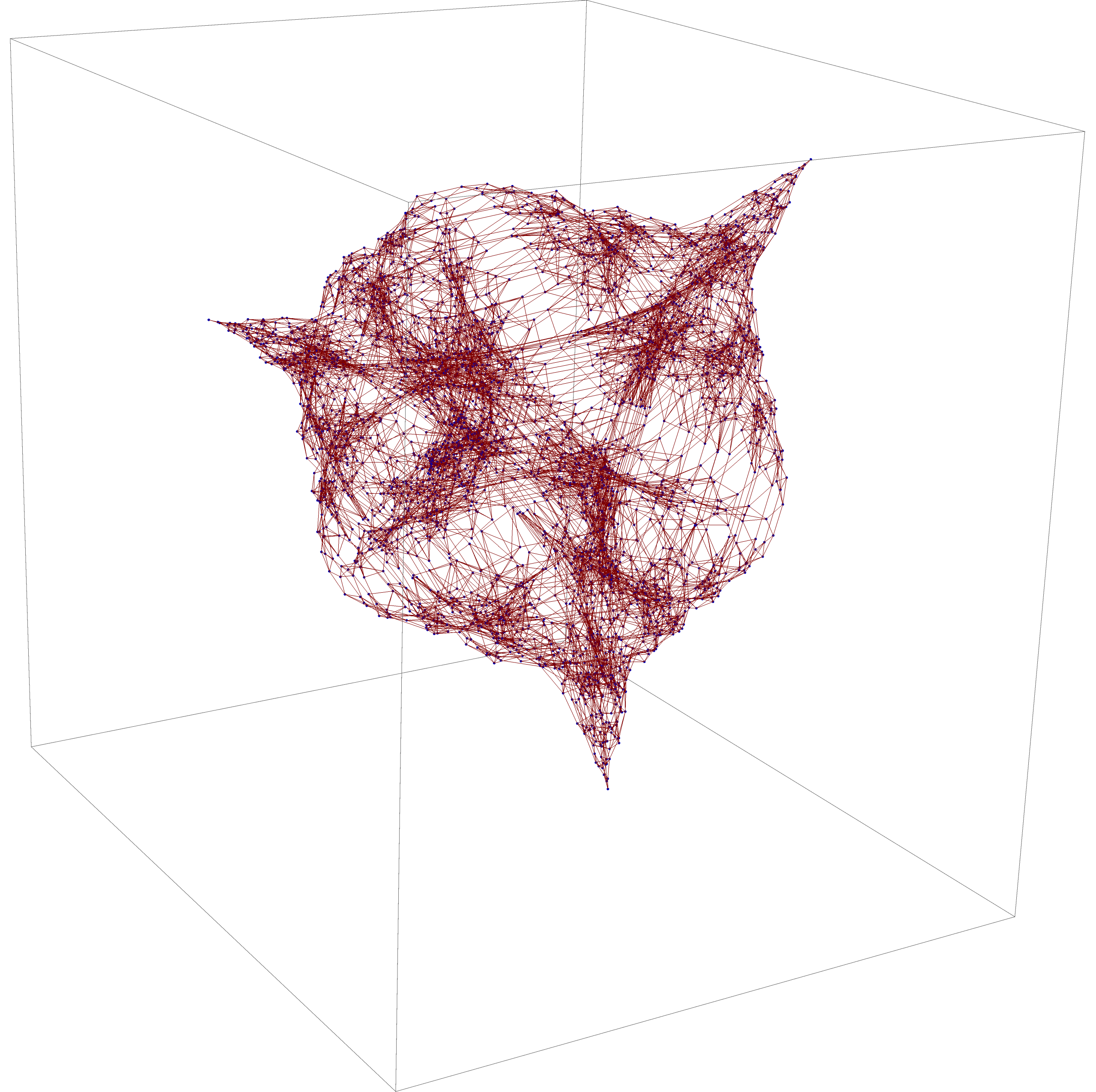}
     \caption{The Schreier graph Hanoi tower group at level 6 }
     \label{Hanoi tower}
\end{figure}
\item[9-] The linear-activity automaton groups is named \textit{Mother group}(Figure \ref{mother group} is used by \cite{bartholdi2010amenability}\cite{amir2013amenability} to prove the amenability of large class of bounded automata groups. The mother group, denoted $\mathfrak{M}_{d,m}$, is defined as the automaton group generated by wreath product \cite{amir2014positive}:
\begin{figure}[h]
\centering
    \includegraphics[width=0.5\textwidth]{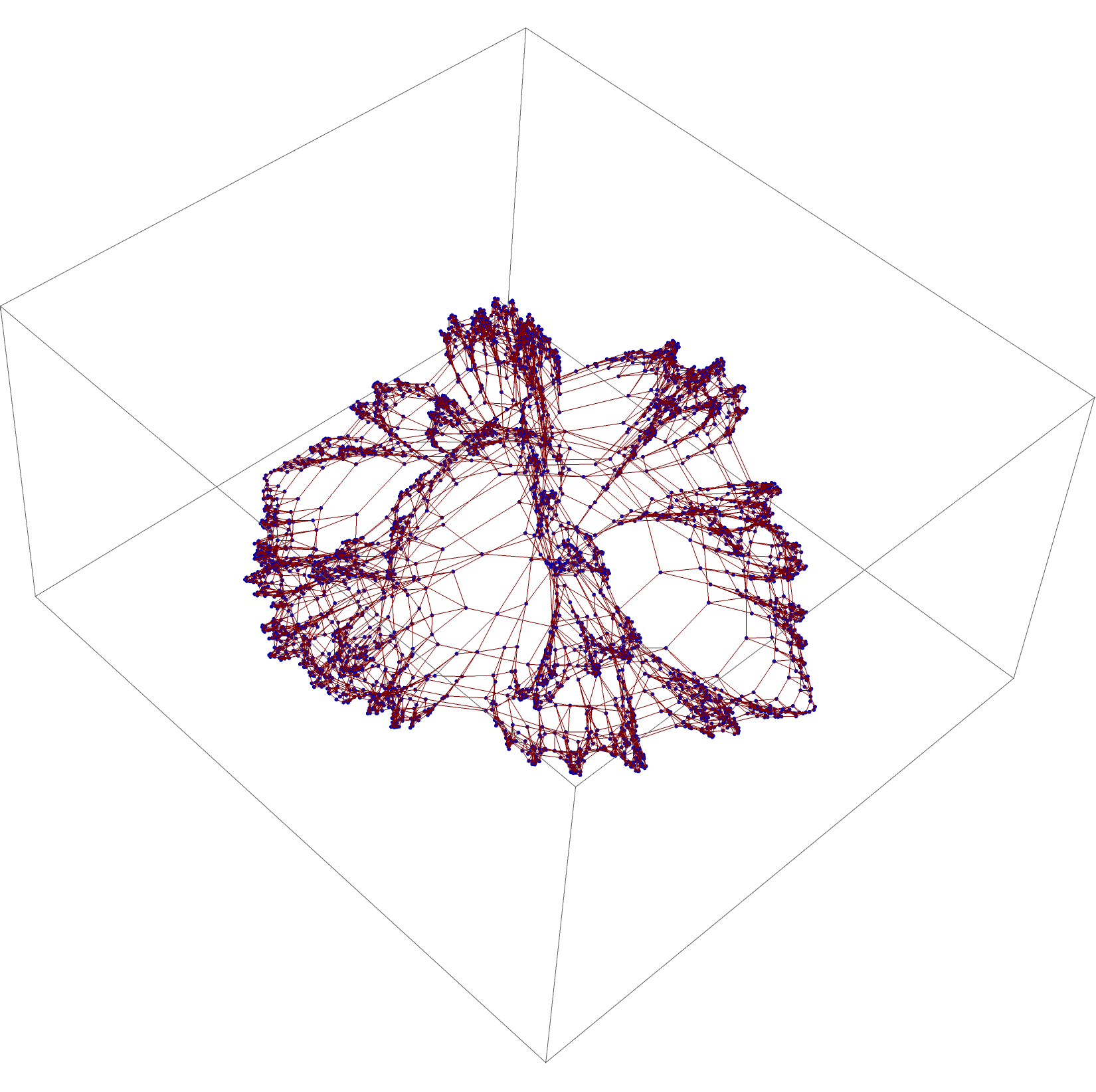}
     \caption{The Schreier graph of the mother group $\mathfrak{M}_{d,m}$ for $d=2,m=2$ at level 12 }
     \label{mother group}
\end{figure}
\begin{flushleft}
$a_{k,\sigma} = \langle a_{k,\sigma}, a_{k-1,\sigma}, 1,\cdots , 1\rangle$  for $0\leq k \leq d$,\\ $a_{-1,\sigma}=\sigma$\\
$b_{k,\rho}=\langle b_{k,\rho}, b_{k-1,\rho}, 1,\cdots , 1\rangle$  for $1\leq k \leq d$,\\ $b_{0,\rho}=\langle b_{0,\rho},1,\cdots,1\rangle$
\end{flushleft}

\end{itemize}
\section{Program and codes}\label{Program codes}

 In this section,  the program code written in Wolfram language by the first author is presented. to compute the adjacency matrixes of Schreier graphs $\Gamma_{n}(G,\Pi, A^{n})$, where $G$ is the group generated by finite closed state automata $\Pi$ and $n\in \mathbb{N}$. We approximate limit $G$-space $\mathcal{J}_{G}$. 
By The comparison of two cases, $G$ be self-similar contracting fractal and self-similar non-contracting fractal groups, in some examples, one can see the vital role of contracting condition influences on limit $G$-space $\mathcal{J}_{G}$ to have a charming asymptotic behavior in geometry, topology, and ergodicity.

The computer program introduced in the Ph.D. thesis \cite{perez2020structural}  written by Perez in 2019 generates Schreier graphs of automatic groups with almost six states and 6-letter of the alphabet. Perez's program gives coincident output with the program proposed in present paper while the above conditions are satisfied. In contrast, the algorithm proposed in this paper has no constraint on the number of states or alphabets. However, the code does its task up to the processing capacity. Indeed this code is a part of a program that builds the Schreier graphs. Furthermore, it computes the spectrum of random walks on Schreier graphs on different levels. Note that this paper is devoted only to the first part of a program that produces Schreier graphs. 
 \subsubsection{automata codes}
 The first part of the program is \textit{Input} section conducted by pice of the code that setups the states of automata and the adjacency structure of the Moor diagram. More percisly 
  The \textit{Input} sector of the program gets:
\begin{itemize}
\item[1 - ] A finite set of alphabet $A$ .
\item[2 - ] The adjacency structure of Moor diagram associated to the automaton. 
\item[4 - ] The permutation action  of  each state $q\in Q$ of automata on alphabet set $A$
\item[3 - ] Integer $n>0$.
\end{itemize}
The \textit{Output} give the adjacency matrix of the Schreier graph at $n$-th level. That is a $\vert A \vert^{n}\times \vert A \vert^{n}$ matrix, such that  $(i,j)$-th entry is an elements of the symbol set $Q\cup \lbrace 0 \rbrace$, where $Q$ is state set of automata $\Pi(Q,A)$.

   for example in following the automata generated by wreath recurrsion:
  \begin{center}
$a=\sigma(b , e)  ,  b=(a , e)  ,  e=(e ,e)$ 
\end{center}

\begin{verbatim}
ClearAll
state[A_, B_] := 
 Function[{u, x}, 
   If[x == 0, Lookup[<|a -> b, b -> a , e -> e|>, u], 
    Lookup[{a -> e, b -> e, e -> e}, u]]][A, B]
G = Graph[{a -> b, b -> a , e -> e, a -> e, b -> e, e -> e}, 
  VertexLabels -> "Name"]
S = VertexList[G]
a[x_] := PermutationReplace[x + 1, Cycles[{{1, 2}}]] - 1
b[x_] := PermutationReplace[x, Cycles[{}]]
e[x_] := PermutationReplace[x, Cycles[{}]]
\end{verbatim}

%%=============================================%%
%% For presentation purpose, we have included  %%
%% \bigskip command. please ignore this.       %%
%%=============================================%%

If the automata have an identity element by use of the following code should remove otherwise, all vertices are adjacent to a trivial edge, $x\mapsto e.x$.

\begin{verbatim}
For[i = 1, i <= Dimensions[S][[1]], i++, If[S[[i]] == e, identity = i]]
identity
S = Delete[S, identity]
\end{verbatim}

Next, generate the vertices of the rooted tree on a given level; in this example, the vertices of the binary rooted tree at level 12 are generated.

\begin{verbatim}

generator = Tuples[{0, 1}, 12];
Clear[u, y, st, H]
H = Table[{}, {i, 1, Dimensions[S][[1]]}];
For[st = 1, st <= Dimensions[S][[1]], st++, 
 For[j = 1; y := {}, j <= Dimensions[generator][[1]], j++, 
  For[i = 1; t = S[[st]]; u = {S[[st]]}, 
   i <= Dimensions[generator][[2]], i++, 
   t = state[t, generator[[j, i]]]; u = Append[u, t]; {i, j}]; 
  y = Append[y, u]]; H[[st]] = y]
\end{verbatim}

The last part is the main part that produces the Schreier graph at a given level. 
\begin{verbatim}
Table[H[[st]][[i, j]][generator[[i, j]]], {st, 1, 
   Dimensions[S][[1]]}, {i, 1, Dimensions[generator][[1]]}, {j, 1, 
   Dimensions[generator][[2]]}];
Schreiergraph = 
 Flatten[Table[
   Table[{StringJoin[
       Table[ToString[generator[[i, j]]], {j, 1, 
         Dimensions[generator][[2]]}]] -> 
      StringJoin[
       Table[ToString[H[[st]][[i, j]][generator[[i, j]]]], {j, 1, 
         Dimensions[generator][[2]]}]], S[[st]]}, {i, 1, 
     Dimensions[generator][[1]]}], {st, 1, Dimensions[S][[1]]}], 1]
GraphPlot[Schreiergraph, DirectedEdges -> False, 
 EdgeRenderingFunction -> Automatic, 
 VertexRenderingFunction -> ({Black, Text[#2, #1]} &)]
Clear[edgespace]
edgespace = {};
For[st = 1, st <= Dimensions[S][[1]], st++,
  For[j = 1, j <= Dimensions[generator][[1]], j++, 
   For[k = 1, k <= Dimensions[generator][[2]], k++, 
    If[H[[st]][[j, k]][generator[[j, k]]] - generator[[j, k]] != 0, 
     edgespace = 
      Append[edgespace, {Table[
         H[[st]][[j, s]][generator[[j, s]]], {s, 1, 
          Dimensions[generator][[2]]}], 
        Table[generator[[j, s]], {s, 1, 
          Dimensions[generator][[2]]}]}]]]]];
edgespace;
ShreierEdge = 
  Union[Table[
    edgespace[[x]][[2]] \[DirectedEdge] edgespace[[x]][[1]], {x, 1, 
     Dimensions[edgespace][[1]]}]];
simplicialSchreieredge = 
  Union[Table[
    StringJoin[
      Table[ToString[edgespace[[x]][[2]][[i]]], {i, 1, 
        Dimensions[generator][[2]]}]] <-> 
     StringJoin[
      Table[ToString[edgespace[[x]][[1]][[i]]], {i, 1, 
        Dimensions[generator][[2]]}]], {x, 1, 
     Dimensions[edgespace][[1]]}]];
gr = Graph[
  Table[Labeled[Schreiergraph[[i, 1]], Schreiergraph[[i, 2]]], {i, 1, 
    Dimensions[Schreiergraph][[1]]}]]
simplicialSchreieredge = SimpleGraph[simplicialSchreieredge]
\end{verbatim}

\bibliographystyle{amsplain}
\bibliography{ref}

\providecommand{\bysame}{\leavevmode\hbox to3em{\hrulefill}\thinspace}
\providecommand{\MR}{\relax\ifhmode\unskip\space\fi MR }
\providecommand{\MRhref}[2]{%
  \href{http://www.ams.org/mathscinet-getitem?mr=#1}{#2}
}
\providecommand{\href}[2]{#2}
\begin{thebibliography}{10}

\bibitem{aleshin1972finite}
Stanislav~V Aleshin, \emph{Finite automata and burnside's problem for periodic
  groups}, Mathematical Notes of the Academy of Sciences of the USSR
  \textbf{11} (1972), no.~3, 199--203.

\bibitem{amir2013amenability}
Gideon Amir, Omer Angel, and B{\'a}lint Vir{\'a}g, \emph{Amenability of
  linear-activity automaton groups}, Journal of the European Mathematical
  Society \textbf{15} (2013), no.~3, 705--730.

\bibitem{amir2014positive}
Gideon Amir and B{\'a}lint Vir{\'a}g, \emph{Positive speed for high-degree
  automaton groups}, Groups, Geometry, and Dynamics \textbf{8} (2014), no.~1,
  23--38.

\bibitem{anantharaman2007entropy}
Nalini Anantharaman, \emph{Entropy and localization of eigenfunctions},
  S{\'e}minaire {\'E}quations aux d{\'e}riv{\'e}es partielles (Polytechnique)
  dit aussi" S{\'e}minaire Goulaouic-Schwartz" (2007), 1--17.

\bibitem{bartholdi2003fractal}
Laurent Bartholdi, Rostislav Grigorchuk, and Volodymyr Nekrashevych, \emph{From
  fractal groups to fractal sets}, Fractals in Graz 2001, Springer, 2003,
  pp.~25--118.

\bibitem{bartholdi2006automata}
Laurent Bartholdi, Andr{\'e}~G Henriques, and Volodymyr~V Nekrashevych,
  \emph{Automata, groups, limit spaces, and tilings}, Journal of Algebra
  \textbf{305} (2006), no.~2, 629--663.

\bibitem{bartholdi2010amenability}
Laurent Bartholdi, Vadim~A Kaimanovich, and Volodymyr~V Nekrashevych, \emph{On
  amenability of automata groups}, Duke Mathematical Journal \textbf{154}
  (2010), no.~3, 575--598.

\bibitem{bartholdi2006thurston}
Laurent Bartholdi and Volodymyr Nekrashevych, \emph{Thurston equivalence of
  topological polynomials}, Acta mathematica \textbf{197} (2006), no.~1, 1--51.

\bibitem{bass2006cyclic}
Hyman Bass, Maria~V Otero-Espinar, Daniel Rockmore, and Charles Tresser,
  \emph{Cyclic renormalization and automorphism groups of rooted trees},
  Springer, 2006.

\bibitem{benjamini2003omega}
Itai Benjamini and Chris Hoffman, \emph{omega-periodic graphs}, arXiv preprint
  math/0308092 (2003).

\bibitem{bondarenko2011ends}
Ievgen Bondarenko, Daniele D’Angeli, and Tatiana Nagnibeda, \emph{Ends of
  schreier graphs of self-similar groups}, preprint (2011).

\bibitem{bondarenko2008classification}
Ievgen Bondarenko, Rostislav Grigorchuk, Rostyslav Kravchenko, Yevgen Muntyan,
  Volodymyr Nekrashevych, Dmytro Savchuk, and Zoran Sunic, \emph{Classification
  of groups generated by 3-state automata over a 2-letter alphabet}, arXiv
  preprint arXiv:0803.3555 (2008).

\bibitem{bourgain2003entropy}
Jean Bourgain and Elon Lindenstrauss, \emph{Entropy of quantum limits},
  Communications in mathematical physics \textbf{233} (2003), no.~1, 153--171.

\bibitem{bridson2013metric}
Martin~R Bridson and Andr{\'e} Haefliger, \emph{Metric spaces of non-positive
  curvature}, vol. 319, Springer Science \& Business Media, 2013.

\bibitem{d2009schreier}
Daniele D'angeli, Alfredo Donno, Michel Matter, and Tatiana Nagnibeda,
  \emph{Schreier graphs of the basilica group}, arXiv preprint arXiv:0911.2915
  (2009).

\bibitem{grigorchuk2014self}
Rostislav Grigorchuk and Dmytro Savchuk, \emph{Self-similar groups acting
  essentially freely on the boundary of the binary rooted tree}, Group theory,
  combinatorics, and computing \textbf{611} (2014), 9--48.

\bibitem{grigorchuk2007schreier}
Rostislav Grigorchuk and Zoran Sunic, \emph{Schreier spectrum of the hanoi
  towers group on three pegs}, arXiv preprint arXiv:0711.0068 (2007).

\bibitem{grigorchuk2006asymptotic}
Rostislav Grigorchuk and Zoran {\v{S}}unik, \emph{Asymptotic aspects of
  schreier graphs and hanoi towers groups}, Comptes Rendus Mathematique
  \textbf{342} (2006), no.~8, 545--550.

\bibitem{grigorchuk2011some}
Rostislav~I Grigorchuk, \emph{Some topics in the dynamics of group actions on
  rooted trees}, Proceedings of the Steklov Institute of Mathematics
  \textbf{273} (2011), no.~1, 64--175.

\bibitem{grigorchuk2002torsion}
Rostislav~I Grigorchuk and Andrzej {\.Z}uk, \emph{On a torsion-free weakly
  branch group defined by a three state automaton}, International Journal of
  Algebra and Computation \textbf{12} (2002), no.~01n02, 223--246.

\bibitem{gromov1981groups}
Michael Gromov, \emph{Groups of polynomial growth and expanding maps (with an
  appendix by jacques tits)}, Publications Math{\'e}matiques de l'IH{\'E}S
  \textbf{53} (1981), 53--78.

\bibitem{gromov1980hyperbolic}
Mikhael Gromov, \emph{Hyperbolic manifolds, groups and actions}, Riemann
  surfaces and related topics: Proceedings of the 1978 Stony Brook Conference
  (State Univ. New York, Stony Brook, NY, 1978), vol.~97, 1980, pp.~183--213.

\bibitem{gromov1987hyperbolic}
\bysame, \emph{Hyperbolic groups}, Essays in group theory, Springer, 1987,
  pp.~75--263.

\bibitem{gromov1992asymptotic}
\bysame, Tech. report, P00001028, 1992.

\bibitem{gromov1999topological}
Misha Gromov, \emph{Topological invariants of dynamical systems and spaces of
  holomorphic maps: I}, Mathematical Physics, Analysis and Geometry \textbf{2}
  (1999), no.~4, 323--415.

\bibitem{jezouin2013quantum}
S{\'e}bastien Jezouin, FD~Parmentier, A~Anthore, U~Gennser, A~Cavanna, Yong
  Jin, and F~Pierre, \emph{Quantum limit of heat flow across a single
  electronic channel}, Science \textbf{342} (2013), no.~6158, 601--604.

\bibitem{kapovich2002boundaries}
Ilya Kapovich and Nadia Benakli, \emph{Boundaries of hyperbolic groups}, arXiv
  preprint math/0202286 (2002).

\bibitem{kigami2001analysis}
Jun Kigami, \emph{Analysis on fractals}, no. 143, Cambridge University Press,
  2001.

\bibitem{lindstrom1990brownian}
Tom Lindstrom and Tom Lindstr{\o}m, \emph{Brownian motion on nested fractals},
  no. 420, American Mathematical Soc., 1990.

\bibitem{AutomGrp1.3.2}
Y.~Muntyan and D.~Savchuk, \emph{{AutomGrp}, automata groups, {V}ersion 1.3.2},
  \href {https://gap-packages.github.io/automgrp}
  {\texttt{https://gap-packages.github.io/}\discretionary
  {}{}{}\texttt{automgrp}}, Sep 2019, Refereed GAP package.

\bibitem{nekrashevych2005self}
Volodymyr Nekrashevych, \emph{Self-similar groups}, no. 117, American
  Mathematical Soc., 2005.

\bibitem{nekrashevych2008groups}
Volodymyr Nekrashevych and Alexander Teplyaev, \emph{Groups and analysis on
  fractals}, Analysis on graphs and its applications \textbf{77} (2008),
  143--180.

\bibitem{nekrashevych2010free}
Volodymyr~V Nekrashevych, \emph{Free subgroups in groups acting on rooted
  trees}, Groups, Geometry, and Dynamics \textbf{4} (2010), no.~4, 847--862.

\bibitem{perez2020structural}
Aitor Perez~Perez, \emph{Structural and spectral properties of schreier graphs
  of spinal groups}, Ph.D. thesis, University of Geneva, 2020.

\bibitem{yue1996ergodic}
Chengbo Yue, \emph{The ergodic theory of discrete isometry groups on manifolds
  of variable negative curvature}, Transactions of the American Mathematical
  Society \textbf{348} (1996), no.~12, 4965--5005.

\end{thebibliography}

\end{document}